\documentclass[12pt]{article}

\usepackage{amssymb,amsfonts,latexsym}
\usepackage{pst-plot}

\newcommand{\case}[4]{\left\{\begin{array}{ll}#1&\mbox{#2}\\#3&\mbox{#4}\end{array}\right.}
\newcommand{\iso}{\cong}
\newcommand{\Cong}{\equiv}
\newcommand{\Mod}{\mathop{\rm mod}\nolimits}
\newcommand{\cF}{{\cal F}}
\newcommand{\cT}{{\cal T}}

\newtheorem{theorem}{Theorem}[section]
\newtheorem{proposition}[theorem]{Proposition}
\newtheorem{lemma}[theorem]{Lemma}

\newcommand{\qed}{  \rule{1ex}{1ex}}
\newenvironment{proof}{\noindent {\bf Proof:}}{{\qed}}

\newcounter{todocounter}
\newcommand{\todo}[1]{
	\addtocounter{todocounter}{1}
	\bigskip
	\noindent{\bf $\ll$ To-do \#\arabic{todocounter}:\rule{10pt}{0pt}#1 $\gg$}\bigskip
}
\newcounter{commentcounter}

\makeatletter
\newfont{\footsc}{cmcsc10 at 8truept}
\newfont{\footbf}{cmbx10 at 8truept}
\newfont{\footrm}{cmr10 at 10truept}
\title{Maximal and Maximum Independent Sets In Graphs With At Most $r$ Cycles}
\author{Bruce E. Sagan\\[-5pt]
\small Department of Mathematics\\[-5pt]
\small Michigan State University\\[-5pt]
\small East Lansing, MI\\[-5pt]
\small \texttt{sagan@math.msu.edu}\\[6pt]
Vincent R. Vatter\thanks{Partially supported by an award from DIMACS and an NSF VIGRE grant to the Rutgers University Department of Mathematics.}\\[-5pt]
\small Department of Mathematics\\[-5pt]
\small Rutgers University\\[-5pt]
\small Piscataway, NJ\\[-5pt]
\small \texttt{vatter@math.rutgers.edu}}

\date{\today \\[6pt]
	\begin{flushleft}
	\small Key Words: cycle, ear decomposition, maximal and maximum independent set\\[6pt]
	\small AMS classification:
	\small Primary 05C35;
	\small Secondary 05C38, 05C69.
	\end{flushleft}
          }

\begin{document}
\maketitle

\begin{abstract}
Let $m(G)$ denote the number of maximal independent sets of vertices
in a graph $G$ and let $c(n,r)$ be the maximum value of $m(G)$ over
all connected graphs with $n$ vertices and at most $r$ cycles.  A
theorem of Griggs, Grinstead, and Guichard gives a formula for
$c(n,r)$ when $r$ is large relative to $n$, while a theorem of Goh,
Koh, Sagan, and Vatter does the same when $r$ is small relative to
$n$.  We complete the determination of $c(n,r)$ for all $n$ and $r$
and characterize the extremal graphs.  Problems for maximum
independent sets are also completely resolved. 
\end{abstract}

\section{Introduction and preliminary lemmas}

Let $G=(V,E)$ be a simple graph.  A subset $I\subseteq V$ is {\em  
independent\/} if there is no edge of $G$ between any two vertices of  
$I$.  Also, $I$ is {\em maximal\/} if it is not properly contained in  
any other independent set.  We let $m(G)$ be the number of maximal  
independent sets of $G$.  Several previous authors have been
interested in the problem of
maximizing $m(G)$ over different families of graphs.

In \cite{mmi1} the authors studied two families of graphs: the family of
all graphs with at most $r$ cycles, and the family of all connected graphs
with at most $r$ cycles.  For the family of all graphs, they were able to
completely settle the problem, by using the result of Moon and
Moser~\cite{mm:cg} (Theorem~\ref{mm} below)
when $n$ is small relative to $r$ and providing
new arguments for all values of $(n,r)$ to which the Moon-Moser Theorem
does not apply (see Theorem~\ref{main1} (I), also below).

For the family of connected graphs, \cite{mmi1}
only characterizes the extremal graphs when $n\ge 3r$
(Theorem~\ref{main1} (II)) while the connected analogue of the
Moon-Moser Theorem (the Griggs-Grinstead-Guichard Theorem, Theorem~\ref{ggg} below)
settles the problem for $n$ small relative to $r$, leaving
a gap between the values where these two theorems apply.
This gap is filled in
Section~\ref{sec-gap} by a careful analysis of the possible endblocks of
extremal graphs.

In the later sections we turn
our attention to maximum independent sets (independent sets of
maximum cardinality) in these two families of graphs.  Like with
maximal independent sets, we start with the case where $n$ is large
relative to $r$ in Section~\ref{sec-maximum} and then consider the
gap in Section~\ref{sec-gap2}.

For the remainder of this section we briefly recount the results we will need.
These results appear in \cite{mmi1} (with the single exception of
Proposition~\ref{endB}, which occurs here in a strengthened form),
and we refer the reader to that paper for examples and proofs.

For any two graphs $G$ and $H$, let $G\uplus H$ denote the disjoint
union of $G$ and $H$, and for any nonnegative integer $t$, let $tG$
stand for the disjoint union of $t$ copies of $G$.  

Let
$$
G(n):=
\left\{\begin{array}{ll}
\frac{n}{3}K_3&\mbox{if $n\Cong0\ (\Mod 3)$,}\\
2K_2\uplus\frac{n-4}{3}K_3&\mbox{if $n\Cong1\ (\Mod
3)$,}\rule{0pt}{20pt}\\
K_2\uplus\frac{n-2}{3}K_3&\mbox{if $n\Cong2\ (\Mod 3)$.}\rule{0pt}{20pt}
\end{array}\right.
$$
Further, let
$$
\mbox{$G'(n):=K_4\uplus\frac{n-4}{3}K_3$ if $n\Cong1\ (\Mod 3)$.}
$$
Also define
$$
g(n):=m(G(n))=
\left\{\begin{array}{ll}
3^\frac{n}{3}&\mbox{if $n\Cong0\ (\Mod 3)$,}\\
4\cdot3^\frac{n-4}{3}&\mbox{if $n\Cong1\ (\Mod 3)$,}\rule{0pt}{20pt}\\
2\cdot3^\frac{n-2}{3}&\mbox{if $n\Cong2\ (\Mod 3)$.}\rule{0pt}{20pt}\\
\end{array}\right.
$$
Note that $m(G'(n))=m(G(n))$ when $n\Cong1\ (\Mod 3)$.

\begin{theorem}[Moon and Moser~\cite{mm:cg}]
\label{mm}
Let $G$ be a graph with $n\ge2$ vertices.  Then
$$
m(G)\le g(n)
$$
with equality if and only if $G\iso G(n)$ or, 
for $n\Cong1\ (\Mod 3)$,
$G\iso G'(n)$.\qed
\end{theorem}

The extremal connected graphs were found by Griggs, Grinstead,  
and Guichard.  To define these graphs we need one more piece of notation.
Let $G$ be a graph all of whose components are complete and let $K_m$
be a complete graph disjoint from $G$.  Construct the graph $K_m * G$ by
picking a vertex $v_0$ in $K_m$ and connecting it to a
single vertex in each component of $G$.  If $n\ge6$ then let
$$
C(n):=
\left\{\begin{array}{ll}
K_3 * \frac{n-3}{3}K_3&\mbox{if $n\Cong0\ (\Mod 3)$,}\\
K_4 * \frac{n-4}{3}K_3&\mbox{if $n\Cong1\ (\Mod 3)$,}\rule{0pt}{20pt}\\
K_4 * \left(K_4\uplus\frac{n-8}{3}K_3\right)
		&\mbox{if $n\Cong2\ (\Mod 3)$.}\rule{0pt}{20pt}
\end{array}\right.
$$
It can be calculated that
$$
c(n):=m(C(n))=
\left\{\begin{array}{ll}
2 \cdot 3^\frac{n-3}{3} + 2^\frac{n-3}{3}&\mbox{if $n\Cong0\ (\Mod
3)$,}\\
3^\frac{n-1}{3} + 2^\frac{n-4}{3}
		&\mbox{if $n\Cong1\ (\Mod 3)$,}\rule{0pt}{20pt}\\
4 \cdot 3^\frac{n-5}{3} + 3 \cdot 2^\frac{n-8}{3}
		&\mbox{if $n\Cong2\ (\Mod 3)$.}\rule{0pt}{20pt}
\end{array}\right.
$$

\begin{theorem}[Griggs, Grinstead, and Guichard~\cite{ggg:nmi}]
\label{ggg}
Let $G$ be a connected graph with $n\ge6$
vertices.  Then
$$
m(G)\le c(n)
$$
with equality if and only if $G\iso C(n)$.\qed
\end{theorem}

The study of $m(G)$ for graphs with a restricted number of cycles began with Wilf.  Let
$$
t(n):=
\case{2^\frac{n-2}{2}+1}{if $n$ is even,}
{2^\frac{n-1}{2}}{if $n$ is odd.}
$$

\begin{theorem}[Wilf~\cite{wil:nmi}]
\label{trees}
If $G$ is a tree with $n\ge 1$ vertices then $m(G)\le t(n)$.\qed
\end{theorem}
Sagan~\cite{sag:nis} gave another proof of this theorem in which he also characterized the extremal graphs, but we will not need them.

Now let
$$
f(n):=2^{\lfloor\frac{n}{2}\rfloor}.
$$
From Theorem~\ref{trees}, one can easily solve the problem for forests.

\begin{theorem}
\label{forests}
If $G$ is a forest with $n\ge 1$ vertices then $m(G)\le f(n)$.\qed
\end{theorem}

To move from trees to a bounded number of cycles, suppose that $n,r$ are positive integers with $n\ge 3r$.  Define
$$
G(n,r):=\case{rK_3\uplus\frac{n-3r}{2}K_2}
{if $n\Cong r\ (\Mod 2)$,}
{(r-1)K_3\uplus\frac{n-3r+3}{2}K_2}
{if $n\not\Cong r\ (\Mod 2)$.\rule{0pt}{20pt}}
$$
Again, it can be computed that
$$
g(n,r) := m(G(n,r)) =
\case{3^r \cdot 2^\frac{n-3r}{2}}
{if $n\Cong r\ (\Mod 2)$,}
{3^{r-1} \cdot 2^\frac{n-3r+3}{2}}
{if $n\not\Cong r\ (\Mod 2)$.\rule{0pt}{30pt}}
$$
It is also convenient to define $G(n,r):=G(n)$ and $g(n,r):=g(n)$ when $n<3r$.  The extremal connected graphs where $n\ge3r$ are given by
$$
C(n,r):=
\left\{\begin{array}{ll}
K_3 * \left((r-1)K_3\uplus\frac{n-3r}{2}K_2\right)
		&\mbox{if $n\Cong r\ (\Mod 2)$,}\\
K_1 * \left(rK_3\uplus\frac{n-3r-1}{2}K_2\right)
		&\mbox{if $n\not\Cong r\ (\Mod 2)$.}\rule{0pt}{20pt}
\end{array}\right.
$$
As usual, we let
$$
c(n,r) := m(C(n,r)) =
\case{3^{r-1} \cdot 2^\frac{n-3r+2}{2} + 2^{r-1}}
{if $n\Cong r\ (\Mod 2)$,}
{3^r \cdot 2^\frac{n-3r-1}{2}}
{if $n\not\Cong r\ (\Mod 2)$.\rule{0pt}{30pt}}
$$

\begin{figure}
\begin{center}
\psset{xunit=0.012in, yunit=0.012in}
\psset{linewidth=1.0\psxunit}
\begin{pspicture}(0,0)(67.5813,97.5452)
\psline(67.5813, 97.5452)(33.7906, 78.0361)     
\psline(67.5813, 58.5271)(33.7906, 78.0361)     
\psline(67.5813, 58.5271)(67.5813, 97.5452)     
\psline(33.7906, 39.0181)(33.7906, 78.0361)     
\psline(0, 58.5271)(33.7906, 78.0361)   
\psline(0, 97.5452)(33.7906, 78.0361)   
\psline(0, 97.5452)(0, 58.5271) 
\psline(33.7906, 0)(33.7906, 39.0181)   
\pscircle*(33.7906, 78.0361){4\psxunit} 
\pscircle*(67.5813, 97.5452){4\psxunit} 
\pscircle*(67.5813, 58.5271){4\psxunit} 
\pscircle*(33.7906, 39.0181){4\psxunit} 
\pscircle*(0, 58.5271){4\psxunit}       
\pscircle*(0, 97.5452){4\psxunit}       
\pscircle*(33.7906, 0){4\psxunit}       
\end{pspicture}
\end{center}
\caption{The exceptional graph $E$}\label{exception}
\end{figure}
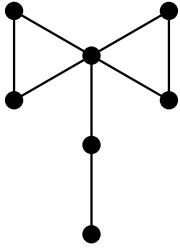

\begin{theorem}[\cite{mmi1}]\label{main1}
Let $G$ be a graph with $n$ vertices and at most $r$ cycles where
$r\ge1$.
\begin{enumerate}
\item[(I)]
If $n\ge 3r-1$ then
$m(G)\le g(n,r)$
with equality if and only if $G\iso G(n,r)$.
\item[(II)]
If $n\ge3r$ then for all such graphs that are connected we have
$m(G)\le c(n,r)$.
Equality occurs if and only if $G\iso C(n,r)$, or if $G$ is one of the
exceptional cases listed in the following table.
$$
\begin{array}{c|c|c}
n	&r	&\mbox{possible $G\not\iso C(n,r)$}\\
\hline
4	&1	& P_4\\
5	&1	& C_5\\
7	&2	& C(7,1), E
\end{array}
$$
(Here $P_4$ and $C_5$ are the path and cycle on 4 and 5 vertices, respectively, and $E$ is the graph shown in Figure~\ref{exception}.)
\end{enumerate}
\end{theorem}

We have a list of inequalities that will be useful in our proofs.
Here and elsewhere it will be convenient to let $g(n,0)=f(n)$ and $c(n,0)=t(n)$.

\begin{lemma}[\cite{mmi1}]\label{inequalities}
We have the following monotonicity results.
\begin{enumerate}
\item[(1)]\label{gnm}
If $r\ge1$ and $n>m\ge3r-1$ then
	$$g(n,r)>g(m,r).$$
\item[(2)]\label{cnm}
If $r\ge1$ and $n>m\ge3r$ then
	$$c(n,r)>c(m,r).$$
\item[(3)]\label{grq}
If $r>q\ge0$ and $n\ge3r-1$ then
	$$g(n,r)\ge g(n,q)$$
with equality if and only if $n$ and $r$ have different parity and
$q=r-1$.
\item[(4)]\label{crq}
If $r>q\ge0$ and $n\ge3r$ then
	$$c(n,r)\ge c(n,q)$$
with equality if and only if $(n,r,q)=(4,1,0)$ or $(7,2,1)$.
\end{enumerate}
\end{lemma}

We also need two results about $m(G)$ for general graphs $G$.  In what  
follows, if $v\in V$ then the {\em open\/} and {\em closed  
neighborhoods of $v$} are $N(v)=\{u\in V\ |\ uv\in E\}$ and
$N[v] = \{v\}\cup N(v)$, respectively.
We also call a block an
{\em endblock of $G$\/} if it has at most one cutvertex in the
graph as a whole.

\begin{proposition}
The invariant $m(G)$ satisfies the following.
\begin{enumerate}
\item[(1)]\label{G-v}
If $v\in V$ then $m(G)\le m(G-v)+m(G-N[v])$.
\item[(2)]\label{endblock}
If $G$ has an endblock $B$ that is isomorphic to a complete graph
then
$$
m(G)=\sum_{v\in V(B)} m(G-N[v]).
$$
\end{enumerate}
In fact, the same equality holds for any complete subgraph $B$
having at least one vertex that is adjacent in $G$ only to other
vertices of $B$.
\end{proposition}
We will refer to the formulas in parts (1) and (2) of this proposition as the
{\it $m$-bound\/} and {\it $m$-recursion\/}, respectively.

Using the fact that the blocks and cutvertices of a graph have a tree structure~\cite[Proposition 3.1.1]{die:gt}, one obtains the following result.

\begin{proposition}
\label{endB}
Every graph has an endblock which intersects at most one non-endblock.  Furthermore, if a graph is not $2$-connected itself, then it contains at least two such endblocks.
\end{proposition}

Note that any block with at least 3 vertices is 2-connected.  Our analysis
of the possible endblocks of the extremal graphs will rely upon Whitney's
Ear Decomposition Theorem from \cite{w:ear}.
\begin{theorem}[Ear Decomposition Theorem]
A graph $B$ is 2-connected if and only if there is a sequence
$$
B_0, B_1,\ldots, B_l=B
$$
such that $B_0$ is a cycle and $B_{i+1}$ is obtained  by
taking a nontrivial path and identifying its two endpoints with two
distinct
vertices of $B_i$. \quad\qed
\end{theorem}
Proofs of the Ear Decomposition Theorem may also be found in Diestel~\cite[Proposition 3.1.2]{die:gt}
and West~\cite[Theorem 4.2.8]{w:gt}.

\section{Filling the gap}\label{sec-gap}

For any graph $G$, it will be convenient to let
$$
r(G)=\mbox{number of cycles of $G$.}
$$

For the family of all graphs on $n$ vertices we have already seen the
maximum value of $m(G)$ for all possible $r(G)$.  Now consider the
family of connected graphs.  If $n\Cong 0\ (\Mod 3)$ then
Theorems~\ref{ggg} and \ref{main1} (II) characterize the maximum for
all possible $r(G)$.  So for the rest of this section we will
concentrate on connected graphs with $n\Cong1,2\ (\Mod 3)$.

Let
$$
r_0:=\lfloor n/3\rfloor=\mbox{the largest value of $r$ for which
Theorem~\ref{main1} (II) is valid.}
$$
Also let
$$
r_1:=r(C(n))=\case{r_0+6}{if $n\Cong1\ (\Mod 3)$,}
{r_0+12}{if $n\Cong2\ (\Mod 3)$.}
$$

To characterize the extremal graphs in the gap $r_0<r<r_1$ we will
need an extension of the star operation.  Let $G$ and $H$
be graphs all of whose components are complete and such that each
component of $H$ has at least 2 vertices.  Construct $K_m * [G,H]$ by
picking a vertex $v_0$ of $K_m$ and connecting it to a single vertex
in each component of $G$ and to two vertices in each component of $H$.
If $n=3r_0+1$ then define
$$
C(n,r):=\case{K_1*[(r_0-1)K_3,K_3]}{if $r=r_0+2$,}
{K_1*[(r_0-2)K_3,2K_3]}{if $r=r_0+4$.}
$$
Note that $C(4,5)$ is not well-defined because then $r_0-2=-1<0$,
and we will leave this graph undefined.  We also need the exceptional graph
$$
C(7,3):=K_1*[3K_2,\emptyset].
$$
For the case $n=3r_0+2$, let
$$
C(n,r):=K_1*[(r_0-1)K_3,2K_2]\quad\mbox{if $r=r_0+1$.}
$$
These will turn out to be the new extremal graphs in the gap.
Two examples may be found in Figure~\ref{Cex2}.

Note that for all graphs just defined we have $m(C(n,r))=m(C(n,r_0))$.
We let
$$
m_0:=m(C(n,r_0))=\case{3^{r_0}}{if $n=3r_0+1$,}
	{4\cdot 3^{r_0-1}+2^{r_0-1}}{if $n=3r_0+2$.}
$$
It will also be convenient to extend the domain of
$c(n,r)$ to all $n$ and $r$ by defining
\begin{eqnarray*}
c(n,r)&=&c(n)\quad\mbox{when $r\ge r(C(n))$,}\\
c(n,r)&=&m_0\quad\mbox{when $r_0<r<r_1$.}
\end{eqnarray*}
We extend the definition of $C(n,r)$ similarly by defining
$C(n,r):=C(n)$ when $r\ge r(C(n))$.

\thicklines
\setlength{\unitlength}{2pt}
\begin{figure}
\begin{center}
\begin{tabular}{ccc}
\psset{xunit=0.012in, yunit=0.012in}
\psset{linewidth=1.0\psxunit}
\begin{pspicture}(0,0)(131.869,128.484)
\psline(26.9165, 33.7906)(65.9346, 33.7906)     
\psline(46.4255, 0)(65.9346, 33.7906)   
\psline(46.4255, 0)(26.9165, 33.7906)   
\psline(85.4436, 0)(65.9346, 33.7906)   
\psline(104.953, 33.7906)(65.9346, 33.7906)     
\psline(104.953, 33.7906)(85.4436, 0)   
\psline(7.40748, 0)(26.9165, 33.7906)   
\psline(7.40748, 0)(46.4255, 0) 
\psline(124.462, 0)(85.4436, 0) 
\psline(124.462, 0)(104.953, 33.7906)   
\psline(98.0785, 128.484)(131.869, 108.975)     
\psline(0, 108.975)(33.7906, 128.484)   
\psline(98.0785, 89.4656)(65.9346, 33.7906)     
\psline(98.0785, 89.4656)(131.869, 108.975)     
\psline(98.0785, 89.4656)(98.0785, 128.484)     
\psline(33.7906, 89.4656)(65.9346, 33.7906)     
\psline(33.7906, 89.4656)(33.7906, 128.484)     
\psline(33.7906, 89.4656)(0, 108.975)   
\pscircle*(65.9346, 33.7906){4\psxunit} 
\pscircle*(26.9165, 33.7906){4\psxunit} 
\pscircle*(46.4255, 0){4\psxunit}       
\pscircle*(85.4436, 0){4\psxunit}       
\pscircle*(104.953, 33.7906){4\psxunit} 
\pscircle*(7.40748, 0){4\psxunit}       
\pscircle*(124.462, 0){4\psxunit}       
\pscircle*(131.869, 108.975){4\psxunit} 
\pscircle*(98.0785, 128.484){4\psxunit} 
\pscircle*(33.7906, 128.484){4\psxunit} 
\pscircle*(0, 108.975){4\psxunit}       
\pscircle*(98.0785, 89.4656){4\psxunit} 
\pscircle*(33.7906, 89.4656){4\psxunit} 
\rput[c](80.9346, 41.7906){$v_0$}
\end{pspicture}
&\rule{30pt}{0pt}&
\psset{xunit=0.012in, yunit=0.012in}
\psset{linewidth=1.0\psxunit}
\begin{pspicture}(0,0)(171.881,133.723)
\psline(47.1363, 31.5663)(85.9406, 35.6448)     
\psline(70.0706, 0)(85.9406, 35.6448)   
\psline(70.0706, 0)(47.1363, 31.5663)   
\psline(101.811, 0)(85.9406, 35.6448)   
\psline(124.745, 31.5663)(85.9406, 35.6448)     
\psline(124.745, 31.5663)(101.811, 0)   
\psline(145.773, 115.77)(171.881, 86.7741)      
\psline(66.4316, 133.723)(105.45, 133.723)      
\psline(0, 86.7741)(26.1082, 115.77)    
\psline(133.716, 78.6618)(85.9406, 35.6448)     
\psline(133.716, 78.6618)(171.881, 86.7741)     
\psline(133.716, 78.6618)(145.773, 115.77)      
\psline(85.9406, 99.9327)(85.9406, 35.6448)     
\psline(85.9406, 99.9327)(105.45, 133.723)      
\psline(85.9406, 99.9327)(66.4316, 133.723)     
\psline(38.1654, 78.6618)(85.9406, 35.6448)     
\psline(38.1654, 78.6618)(26.1082, 115.77)      
\psline(38.1654, 78.6618)(0, 86.7741)   
\pscircle*(85.9406, 35.6448){4\psxunit} 
\pscircle*(47.1363, 31.5663){4\psxunit} 
\pscircle*(70.0706, 0){4\psxunit}       
\pscircle*(101.811, 0){4\psxunit}       
\pscircle*(124.745, 31.5663){4\psxunit} 
\pscircle*(171.881, 86.7741){4\psxunit} 
\pscircle*(145.773, 115.77){4\psxunit}  
\pscircle*(105.45, 133.723){4\psxunit}  
\pscircle*(66.4316, 133.723){4\psxunit} 
\pscircle*(26.1082, 115.77){4\psxunit}  
\pscircle*(0, 86.7741){4\psxunit}       
\pscircle*(133.716, 78.6618){4\psxunit} 
\pscircle*(85.9406, 99.9327){4\psxunit} 
\pscircle*(38.1654, 78.6618){4\psxunit} 
\rput[c](102.9406, 40.6448){$v_0$}
\end{pspicture}
\\\\
C(13,8)
&&
C(14,5)
\end{tabular}
\caption{Examples of $C(n,r)$ for $n< 3r$}\label{Cex2}
\end{center}
\end{figure}
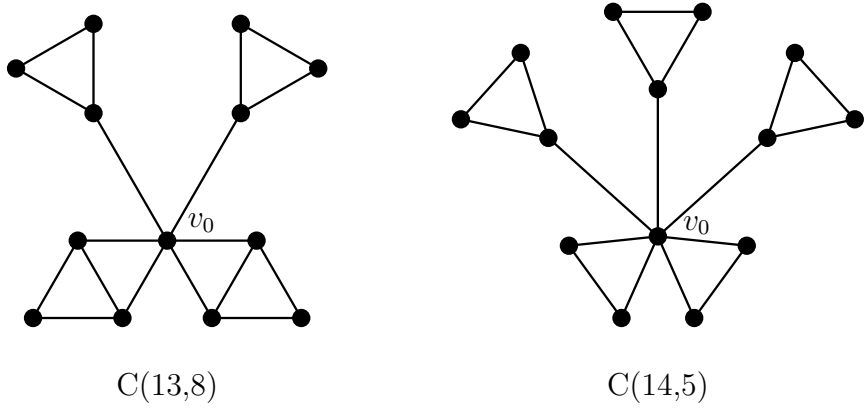

We will need a special case of the Moon-Moser
transformation~\cite{ggg:nmi,mm:cg} which, in conjunction with
the Ear Decomposition Theorem,
will prove useful in cutting down on the
number of cases to consider.
Suppose the graph $G$ contains a path $tuvw$ such that
\begin{enumerate}
\item[(a)] $\deg u=\deg v=2$,
\item[(b)] the edge $uv$ lies on a cycle, and
\item[(c)] $u$ is in at least as many maximal independent sets as $v$.
\end{enumerate}
Then construct the (connected) graph $G_{u,v}$ where $V(G_{u,v})=V(G)$ and
$$
E(G_{u,v})=E(G)\cup\{tv\}-\{vw\}.
$$
The edge $uv$ lies on a unique cycle (a $3$-cycle) in $G_{u,v}$, and so
$r(G_{u,v})\le r(G)$ by (b).

\begin{lemma}\label{2deg2}
Suppose $G$ contains a path $tuvw$ satisfying (a)--(c).  Then
$m(G_{u,v})\ge m(G)$,
with equality only if $N_G(w)-\{v\}\subseteq N_G(t)$.
\end{lemma}
\begin{proof}
Let $M$ be a maximal independent set in $G$.  Then there are three
mutually exclusive possibilities for $M$, namely $u,v\not\in M$;
$u\in M$ and $v\not\in M$; or $v\in M$ and $u\not\in M$.  In the first
two cases, $M$ gives rise to distinct maximal independent set(s) in
$G_{u,v}$
as in the following chart.
$$
\begin{array}{c|c}
\mbox{type of MIS $M$ in $G$}&\mbox{corresponding
MIS(s) in $G_{u,v}$}\\
\hline\hline
u,v\not\in M		&	M\\
\hline
u\in M, v\not\in M	&	M, M\cup\{v\}-\{u\}\\
\end{array}
$$

Therefore by (c) we have $m(G_{u,v})\ge m(G)$ without even considering
the
third case where $v\in M$ and $u\not\in M$.  In this case, if there is
some
vertex $x\in N_G(w)-N_G(t)-\{v\}$, then there is a maximal independent
set $M$ in $G$ with $\{t,x,v\}\subseteq M$.  Hence $M-\{v\}$ is a
maximal independent set in $G_{u,v}$ so $m(G_{u,v})>m(G)$, as desired.
If such
a vertex $x$ does not exist, we have $N_G(w)-\{v\}\subseteq
N_G(t)$.
\end{proof}

Given a graph $G$, we let $\mathcal{T}(G)$ denote the set of all graphs that can be obtained from $G$ by applying a maximal sequence of these special Moon-Moser transformations.  By Lemma~\ref{2deg2}, every graph in $\mathcal{T}(G)$ has at least as many cycles and at least as many maximal independent sets as $G$.  Furthermore, since all of these graphs are formed by maximal sequences of transformations, if $H\in\mathcal{T}(G)$ then $H$ cannot contain a path $tuvw$ satisfying (a)--(c) above.  Another way to state this is that every $H\in\mathcal{T}(G)$ has the following property:
\begin{enumerate}
\item[($\Delta$)] If $uv\in E(H)$ lies on a cycle and $\deg u=\deg
v=2$, then $uv$ lies on a $3$-cycle.
\end{enumerate}

Before closing the gap, we wish to
mention a result which we will need to rule out some graphs
from the list of possible extremals.  To state this lemma,
we say that a vertex $v\in V(G)$ is {\it duplicated\/} if there is a
vertex $w\in V(G)$ such that $v$ and $w$
have the same neighbors, that is, $N(v)=N(w)$.
\begin{lemma}
\label{dup}
Let $G$ be a graph with $n$ vertices and a vertex $v$ that is duplicated.
\begin{enumerate}
\item[(1)]  We have $m(G)=m(G-v)$.
\item[(2)]  If $n\ge 7$, $n\Cong 1\mbox{ or } 2\ (\Mod 3)$, and $G$
is connected with less than $r_1$ cycles then
$m(G)<m_0$.
\end{enumerate}
\end{lemma}
\begin{proof}
If $u$ and $v$
are duplicated vertices, then they lie in the same maximal
independent sets and neither is a cutvertex.   So $m(G)=m(G-v)$ and
under the hypotheses of (2), $m(G-v)\le c(n-1)<m_0$.
\end{proof}

We now finish our characterization of the extremal graphs.
\begin{theorem}\label{gap}
Let $G$ be a connected graph with $n\ge 7$ vertices,
$n\Cong 1\mbox{ or } 2\ (\Mod 3)$, and less than $r_1$ cycles.
Then
$$
m(G)\le m_0,
$$
with equality if and only if $G\iso C(n,s)$ for some $s$ with
$r_0\le s<r_1$.
\end{theorem}

\begin{proof}
We prove the theorem by induction on $n$.  The cases where $n\le 10$
have been checked by computer, so let $G$ be an extremal connected graph with 
$n>10$ vertices and less than $r_1$ cycles.

First note that it suffices to prove the theorem for graphs that satisfy 
($\Delta$):  If $G$ satisfies the hypotheses of the theorem then every 
graph $H\in\mathcal{T}(G)$ also satisfies these hypotheses and
satisfies $(\Delta)$.  Since $G$ is extremal, then in
Lemma~\ref{2deg2} we would always have equality and thus the given subset
relation,  when replacing $G$ by $G_{u,v}$.  Since none of our
candidate extremal graphs can be generated by this transformation
if such a condition is imposed, we must have that $\cT(G)=\{G\}$ and
so $G$ satisfies ($\Delta$).

Pick an endblock $B$ of $G$ satisfying the conclusion of 
Proposition~\ref{endB} with $|V(B)|$ maximum among all such endblocks.  If 
$G\neq B$ then we will use $x$ to denote the cutvertex of $G$ in $B$.  
The argument depends on the nature of $B$.

If $B\iso K_2$ then the argument used in the proof of
Theorem~\ref{main1} (II) can be easily adapted for use in this
context.  Let $V(B)=\{x,v\}$  so that
$\deg v=1$ and $\deg x\ge2$.  By the choice of $B$,
$G-N[v]$ is the union of some number of  
$K_1$'s and a connected graph with at most $n-2$ vertices and at most  
$r$ cycles.  Also, $G-N[x]$ has at most $n-3$ vertices and at most $r$  
cycles, so the $m$-recursion and monotonicity give, for $n\ge 11$,
$$
m(G)\le c(n-2,r) + g(n-3,r)\le c(n-2)+g(n-3)\le m_0,
$$
with equality if and only if $n=3r_0+2$ and $G\iso C(n,r_0)$.

All other possible endblocks must be 2-connected and so we will
use the Ear Decomposition Theorem to
organize the cases to consider based on $l$, the number of paths that
are added to the initial cycle.

If $l=0$ then ($\Delta$) guarantees that $B\iso K_3$.  Let
$V(B)=\{x,v,w\}$ where $x$ denotes the cutvertex.  Let $i$ denote the
number of other $K_3$ endblocks containing $x$.  If the graph consists
entirely of $K_3$ endblocks which intersect at $x$, then 
$$
m(G)=2^{\frac{n-1}{2}}+1,
$$
which shows that $G$ is not extremal.  Thus $x$ is adjacent to at
least one vertex which does not lie in a $K_3$ endblock.  Since $B$
was chosen with $|V(B)|$ maximal, it follows that $G-N[v]=G-N[w]$ has
some number of trivial components, $i$ components isomorphic to $K_2$,
and at most one other component, $H$, with at most $n-2i-3$ vertices
at at most $r-i-1\le r_1-i-2$ cycles.  Since $x$ is adjacent to at least one
vertex not in the $K_3$ endblocks, $G-N[x]$ has at most $n-2i-4$
vertices.  This gives us the upper bound 
\begin{equation}
\label{K3}
m(G)\le 2^{i+1}c(n-2i-3,r_1-i-2)+g(n-2i-4).
\end{equation}

To show that this bound is always at most $m_0$, we consider the 
two values $n=3r_0+1$ and $n=3r_0+2$ as well as the three
possible congruence classes of $i$ modulo 3 separately.  So let 
$j=\lfloor i/3\rfloor$.  Considering the number of vertices in
$G-N[x]$ gives $n-2i-4\ge0$ and translating this into a bound
involving $r_0$ and $j$ gives $r_0\ge 2j+k_0$ where 
$1\le k_0\le 3$ depending on which of the six cases we are in.
We now wish to show that the right-hand side
of~(\ref{K3}) is a strictly decreasing function of $j$ for any fixed
but sufficiently large $n$.
This is clearly true of the $g(n-2i-4)$ term, so let 
$f(r_0,j)=2^{i+1} c(n-2i-3,r_1-i-2)$ where the right side has been
converted to a function of $r_0$ and $j$.  In all cases, we get that
$$
f(r_0,j)=a(r_0)2^j+b(r_0) (8/9)^j
$$
for certain functions $a(r_0),b(r_0)$.  It follows that 
$f(r_0,j)-f(r_0,j+1)>0$ if and only if $b(r_0)/a(r_0)>9(3/2)^{2j}$.
Solving for $r_0$ shows that we have a decreasing function of $j$ for
$r_0\ge 2j+k_1$ where $4\le k_1\le 7$.  So it suffices to check that
the right-hand side of~(\ref{K3}) is at most $m_0$ for 
$2j+k_0\le r_0 \le 2j +k_1$.  This is done by substituting each value
of $r_0$ in turn to get a function of $j$ alone, noting that this
function is decreasing for all $j$ sufficiently large to make
$r_0\ge3$, and then verifying that this function is bounded by $m_0$
when $j$ is at this minimum value.  The only cases where we get equality
are when $n=3r_0+2$ and $G\iso C(n,r_0)$.

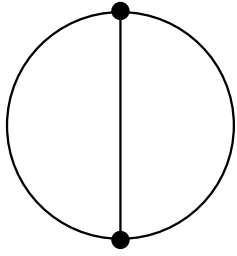
\begin{figure}
\begin{center}
\psset{xunit=0.012in, yunit=0.012in}
\psset{linewidth=1.0\psxunit}
\begin{pspicture}(-50,-50)(50,50)
\pscircle(0,0){50\psxunit}      
\psline(0, -50)(0, 50)        
\pscircle*(0, 50){4\psxunit}   
\pscircle*(0, -50){4\psxunit}  
\end{pspicture}
\caption{$D$, the multigraph for $l=1$}\label{D}
\end{center}
\end{figure}

If $l=1$, then $B$ must be a subdivision of the multigraph $D$ in
Figure~\ref{D}, i.e., it must be obtained from $D$ by inserting
vertices of degree 2 into the edges of $D$.  By ($\Delta$),
we can insert at most one vertex into an edge, unless one of the
inserted vertices is the cutvertex $x$ in which case it is possible
to insert a vertex before and after $x$ as well.  To turn
this multigraph into a graph, it is necessary to subdivide at least
two of the edges.  If all three edges are subdivided, or two edges are
subdivided and $x$ is one of the original vertices of $D$, then $G$
has a duplicated vertex and so is not extremal by Lemma~\ref{dup}.
In the only remaining case, the following lemma applies.

\begin{lemma}[Triangle Lemma]
\label{tri}
Suppose $G$ contains three vertices $\{u,v,w\}$ satisfying the
following restrictions.
\begin{enumerate}
\item[(a)] These vertices form a $K_3$ with $\deg u=2$ and
$\deg v,\deg w\ge 3$.
\item[(b)] The graph $G-\{u,v,w\}$ is connected.
\end{enumerate}
Then
$m(G)\le m_0$
with equality only if $n=3 r_0+1$ and $G\iso C(n,s)$ for some $s$ with
$r_0\le s<r_1$.
\end{lemma}
\begin{proof}
Because of~(a), the $K_3$ satisfies the
alternative hypothesis in the $m$-recursion.  Using induction
to evaluate the $c$ and $g$ functions, we get
\begin{eqnarray*}
m(G)&=&m(G-N[u])+m(G-N[v])+m(G-N[w])\\
	&\le&c(n-3,r-1)+2g(n-4,r-1)\\
	&\le&m_0.
\end{eqnarray*}
with equality only if $n=3r_0+1$, $G-N[v]\iso G-N[w]\iso G(n-4)$,
and $G-N[u]\iso C(n-3,s)$ for some $s\le r_0+4$.  These
easily imply the conclusion of the lemma.
\end{proof}

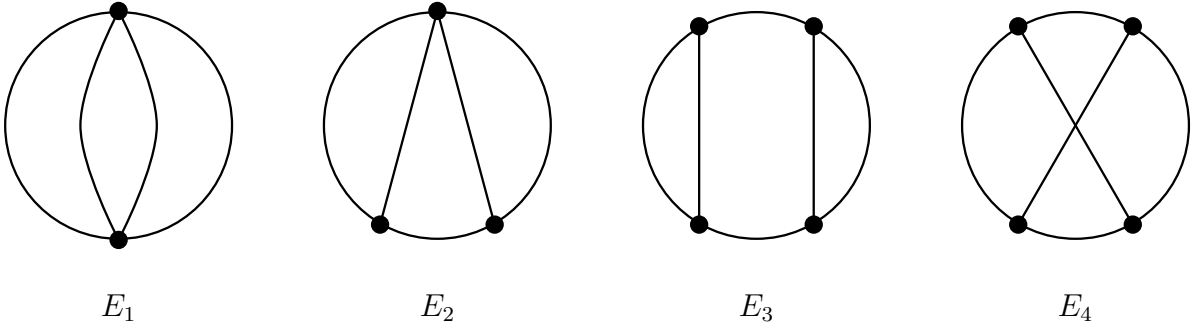
\begin{figure}
\begin{tabular}{ccccccc}
\psset{xunit=0.012in, yunit=0.012in}
\psset{linewidth=1.0\psxunit}
\begin{pspicture}(-50,-50)(50,50)
\pscircle(0,0){50\psxunit}      
\pscurve(0,50)(16.6667,0)(0,-50)
\pscurve(0,50)(-16.6667,0)(0,-50)
\pscircle*(0, 50){4\psxunit}   
\pscircle*(0, -50){4\psxunit}  
\end{pspicture}
&\rule{10pt}{0pt}&
\psset{xunit=0.012in, yunit=0.012in}
\psset{linewidth=1.0\psxunit}
\begin{pspicture}(-50,-50)(50,50)
\pscircle(0,0){50\psxunit}      
\psline(25, -43.3013)(0, 50)    
\psline(-25, -43.3013)(0, 50)   
\pscircle*(0, 50){4\psxunit}    
\pscircle*(25, -43.3013){4\psxunit}     
\pscircle*(-25, -43.3013){4\psxunit}    
\end{pspicture}
&\rule{10pt}{0pt}&
\psset{xunit=0.012in, yunit=0.012in}
\psset{linewidth=1.0\psxunit}
\begin{pspicture}(-50,-50)(50,50)
\pscircle(0,0){50\psxunit}      
\psline(25, -43.3013)(25, 43.3013)      
\psline(-25, -43.3013)(-25, 43.3013)    
\pscircle*(25, 43.3013){4\psxunit}      
\pscircle*(-25, 43.3013){4\psxunit}     
\pscircle*(25, -43.3013){4\psxunit}     
\pscircle*(-25, -43.3013){4\psxunit}    
\end{pspicture}
&\rule{10pt}{0pt}&
\psset{xunit=0.012in, yunit=0.012in}
\psset{linewidth=1.0\psxunit}
\begin{pspicture}(-50,-50)(50,50)
\pscircle(0,0){50\psxunit}      
\psline(25, -43.3013)(-25, 43.3013)     
\psline(-25, -43.3013)(25, 43.3013)     
\pscircle*(25, 43.3013){4\psxunit}      
\pscircle*(-25, 43.3013){4\psxunit}     
\pscircle*(25, -43.3013){4\psxunit}     
\pscircle*(-25, -43.3013){4\psxunit}    
\end{pspicture}
\\\\
$E_1$&&$E_2$&&$E_3$&&$E_4$
\end{tabular}
\caption{The multigraphs for $l=2$}\label{E}
\end{figure}

For $l\ge2$, we must consider the two congruence classes for $n$
separately.  First consider $n=3r_0+1$.   The following lemma will help
eliminate many cases.  In it, we use $r(v)$ to denote the number of
cycles of $G$ containing the vertex $v$.

\begin{lemma}
\label{36}
Let $n=3r_0+1$ and let $v$ be a non-cutvertex with $\deg v\ge3$ and
$r(v)\ge6$.  Then $G$ is not extremal.
\end{lemma}
\begin{proof}
Using the $m$-bound, we have
\begin{eqnarray*}
m(G)&\le&c(n-1,r-6)+g(n-4,r-6)\\
	&\le&c(3r_0,r_0-1)+g(3r_0-3,r_0-1)\\
	&=&2\cdot 3^{r_0-1}+3^{r_0-1}\\
	&=&m_0.
\end{eqnarray*}
However, if we have equality then this forces us to have $G-v\iso
C(3r_0,r_0-1)$ and $G-N[v]\iso G(3r_0-3,r_0-1)$.  The only way this can 
happen is if
$G\iso C(3r_0+1,r_0+2)$ where $v$ is one of the degree 3 vertices in
the 4-vertex block.  But then $v$ is in only 3 cycles, contradicting
our hypothesis that $r(v)\ge6$.
\end{proof}

When $l=2$, $B$ must be a subdivision of one of the
multigraphs in Figure~\ref{E}.  Lemma~\ref{36} shows that $B$ cannot be a subdivision of $E_1$, $E_4$, or any block formed by a sequence of length $l\ge3$ since in all these cases there are at least two vertices having degree at least 3 and lying in at least 6 cycles.  So even if $B$ has a cutvertex, there will 
still be a non-cutvertex in $B$ satisfying the hypotheses of the lemma.

If $B$ is formed by subdividing $E_2$, the same lemma shows that we need
only consider the case where the vertex of degree $4$ in $E_2$ is a cutvertex, $x$, of $G$.
Also, since $B$ can't
have duplicated vertices and must satisfy ($\Delta$), each pair of doubled edges has a vertex
inserted
in exactly one edge.  This means there are only two
possibilities for $B$, depending on whether the non-doubled edge is
subdivided or not, and it is easy to check that in both cases $G$ is
not extremal by using the $m$-bound on the vertex $x$.

Finally, if $B$ is a subdivision of $E_3$ then, because of the pair of 
disjoint doubled
edges, there will always be one doubled edge which does not contain a cutvertex.  In $B$ that pair will give rise to either a duplicated vertex
or a $K_3$ satisfying the hypotheses of the Triangle Lemma, and thus in 
either
case we will be done.  This ends the proof for $n=3r_0+1$.

%
%

Now we look at the case where $n=3r_0+2$.  The analogue of
Lemma~\ref{36} in this setting is as follows and since the proof is
similar, we omit it.

\begin{lemma}
\label{312}
Let $n=3r_0+2$ and let $v$ be a non-cutvertex that satisfies either
\begin{enumerate}
\item[(1)] $\deg v\ge 3$ and $r(v)\ge 12$, or
\item[(2)] $\deg v\ge 4$ and $r(v)\ge 6$.
\end{enumerate}
Then $G$ is not extremal.
\end{lemma}

The ideas used to rule out $E_3$ for $n=3r_0+1$ will be
used many times in the current case, so we codify them in the lemma 
below.

\begin{lemma}
\label{pairs}
Suppose $n=3r_0+2$ and the block $B$ is a subdivision of a multigraph
having two
disjoint submultigraphs  each of which is of one of the following forms:
\begin{enumerate}
\item[(i)] a doubled edge, or
\item[(ii)] a vertex $v$ satisfying $\deg v\ge 3$ and $r(v)\ge 12$, or
\item[(iii)] a vertex $v$ satisfying $\deg v\ge 4$ and $r(v)\ge 6$.
\end{enumerate}
Then $G$ is not extremal.
\end{lemma}
\begin{proof}
If any set of doubled edges has both edges subdivided exactly once,
then $G$ is not extremal by Lemma~\ref{dup}.
Otherwise, since $B$ has at most one cutvertex $x$ in $G$, either the
hypotheses of
the Triangle Lemma or of the previous lemma will be satisfied.
\end{proof}

Finally, we will need a way to eliminate blocks that only have
vertices of degree at most 3, but not sufficiently many cycles to
satisfy Lemma~\ref{312}~(1).  One way would be to make sure that
$G-N[v]$ is connected.  Since a given multigraph $M$ has many possible
subdivisions, we also need a criterion on $M$ that
will guarantee that most of the subdivisions will have the desired
connectivity.

\begin{lemma}
\label{con}
Let $n=3r_0+2$.
\begin{enumerate}
\item[(1)] Suppose that $G$ contains a non-cutvertex $v$ such that
$\deg v\ge3$ and $r(v)\ge6$.
Suppose further that $G-N[v]$ contains at most two nontrivial
components and that if there are two, then one of them is a star (a complete bipartite graph of the form $K_{1,s}$).  Then
$G$ is not extremal.
\item[(2)]  Suppose $G$ comes from subdivision of a multigraph $M$
that contains a vertex $v$ with $\deg v=3$, $r(v)\ge 6$, and such
that all vertices in $N_M[v]$
are non-cutvertices in $M$ and $M-N_M[v]$ is connected.
Suppose further that there are at most two edges of $M$ between the
elements of $N_M(v)$.  Then $G$ and $v$ satisfy the hypotheses of~(1).
\end{enumerate}
\end{lemma}
\begin{proof}
For (1), first assume that there is only one nontrivial component in
$G-N[v]$.
Then, using the induction hypothesis about the behavior of graphs in
the gap,
\begin{eqnarray*}
m(G)	&\le& c(n-1,r-6)+c(n-4,r-6)\\
	&=& c(3r_0+1,r_0+5)+c(3r_0-2,r_0+5)\\
	&=& c(3r_0+1,r_0)+c(3r_0-2)\\
	&=& 3^{r_0} + 3^{r_0-1} + 2^{r_0-2}\\
	&<& m_0.
\end{eqnarray*}

In the case with a star component, we use
$m(G)\le c(n-1,r-6)+2c(n-6,r-6)$ to obtain the same result.

For (2) we will break the proof into several cases depending on how
the edges of $M$ at $v$ are subdivided in $G$, noting that by our hypotheses
each can be subdivided at most once and that the same is true of any edge between elements of $N_M(v)$.  Let $N_M(v)=\{s,t,u\}$ and let $H$ be the subdivision of $L=M-N_M[v]$ induced by $G$.  Note that
$H$ is connected by assumption.  If none of $vs,vt,vu$ are subdivided
in $G$, then $G-N_G[v]$ is just $H$ together, possibly, with some vertices of degree one attached (if any edges from $s$, $t$, or $u$ to $L$ were subdivided) and some trivial components (if any edges between $s$, $t$, and $u$ were subdivided).
If exactly one of the three edges is subdivided,
suppose it is $vs$.  Then it is possible that $s$ is in a different
nontrivial component of $G-N_G[v]$ than $H$.  But since there are at most two edges
from $s$ to $t$ and $u$, the component of $s$ is a star.  Now suppose that $vs$ and $vt$ are subdivided, but not
$vu$.  At
least one of $s,t$ are connected to $H$ otherwise $u$ becomes a
cutvertex.  So again the only possibility for a nontrivial component
other than $H$ is a star containing either $s$ or $t$, but not both.
Finally, if all three edges are subdivided, then $G-N_G[v]$ is
connected because $v$ is not a cutvertex in $M$.
\end{proof}

We need a little terminology before we handle the $n=3r_0+2$
case.  Let $L$ and $M$ be 2-connected multigraphs with no vertices of
degree
two.  We say that {\it $M$ is a child of $L$\/} if there is some
sequence $B_0,B_1,\ldots,B_l$ formed as in the Ear Decomposition Theorem with
$B_{l-1}=L$ and $B_l=M$.  We will use words like ``descendant,''
``parent,'' and so on in a similar manner.

We now pick up the proof for $n=3r_0+2$ where we left off, namely with
$l=2$.  Lemma~\ref{pairs}~(iii) shows that $B$ cannot be a subdivision
of $E_1$ or any of its
descendants.  Also, if $B$ is a subdivision of $E_2$, then by
Lemma~\ref{312}~(2) we need only consider the case where the vertex of
degree 4 is a
cutvertex $x$, and the same argument we used in the $3r_0+1$ case shows
that such graphs are not extremal.

\thicklines
\begin{figure}
\begin{center}
\begin{tabular}{ccccccc}
\psset{xunit=0.012in, yunit=0.012in}
\psset{linewidth=1.0\psxunit}
\begin{pspicture}(-50,-64)(50,64)
\pscircle(0,0){50\psxunit}      
\psline(25, -43.3013)(0, 50)    
\psline(-25, -43.3013)(0, 50)   
\psline(25, 43.3013)(0, -50)    
\pscircle*(0, 50){4\psxunit}    
\pscircle*(25, -43.3013){4\psxunit}     
\pscircle*(-25, -43.3013){4\psxunit}    
\pscircle*(0, -50){4\psxunit}   
\pscircle*(25, 43.3013){4\psxunit}      
\rput[c](0,64){$x$}
\rput[c](0,-64){$v$}
\end{pspicture}
&\rule{10pt}{0pt}&
\psset{xunit=0.012in, yunit=0.012in}
\psset{linewidth=1.0\psxunit}
\begin{pspicture}(-50,-64)(50,64)
\pscircle(0,0){50\psxunit}      
\psline(25, -43.3013)(0, 50)    
\psline(-25, -43.3013)(0, 50)   
\psline(50, 0)(-50, 0)  
\pscircle*(0, 50){4\psxunit}    
\pscircle*(25, -43.3013){4\psxunit}     
\pscircle*(-25, -43.3013){4\psxunit}    
\pscircle*(-50, 0){4\psxunit}   
\pscircle*(50, 0){4\psxunit}    
\rput[c](0,64){$x$}
\end{pspicture}
&\rule{10pt}{0pt}&
\psset{xunit=0.012in, yunit=0.012in}
\psset{linewidth=1.0\psxunit}
\begin{pspicture}(-50,-64)(50,64)
\pscircle(0,0){50\psxunit}      
\psline(-35.3553, -35.3553)(0, 50)      
\psline(0, -50)(0, 50)  
\psline(35.3553, -35.3553)(0, 50)       
\pscircle*(0, 50){4\psxunit}    
\pscircle*(-35.3553, -35.3553){4\psxunit}       
\pscircle*(0, -50){4\psxunit}   
\pscircle*(35.3553, -35.3553){4\psxunit}        
\rput[c](0,64){$x$}
\end{pspicture}
&\rule{10pt}{0pt}&
\psset{xunit=0.012in, yunit=0.012in}
\psset{linewidth=1.0\psxunit}
\begin{pspicture}(-50,-64)(50,64)
\pscircle(0,0){50\psxunit}      
\psline(-40.4508, -29.3893)(0, 50)      
\psline(-15.4508, -47.5528)(0, 50)      
\psline(15.4508, -47.5528)(0, 50)       
\psline(40.4508, -29.3893)(0, 50)       
\pscircle*(0, 50){4\psxunit}    
\pscircle*(-40.4508, -29.3893){4\psxunit}       
\pscircle*(-15.4508, -47.5528){4\psxunit}       
\pscircle*(15.4508, -47.5528){4\psxunit}        
\pscircle*(40.4508, -29.3893){4\psxunit}        
\rput[c](0,64){$x$}
\end{pspicture}
\\\\
$F_1$&&$F_2$&&$F_3$&&$G_1$
\end{tabular}
\caption{Three children and one grandchild of $E_2$}\label{E2d}
\end{center}
\end{figure}
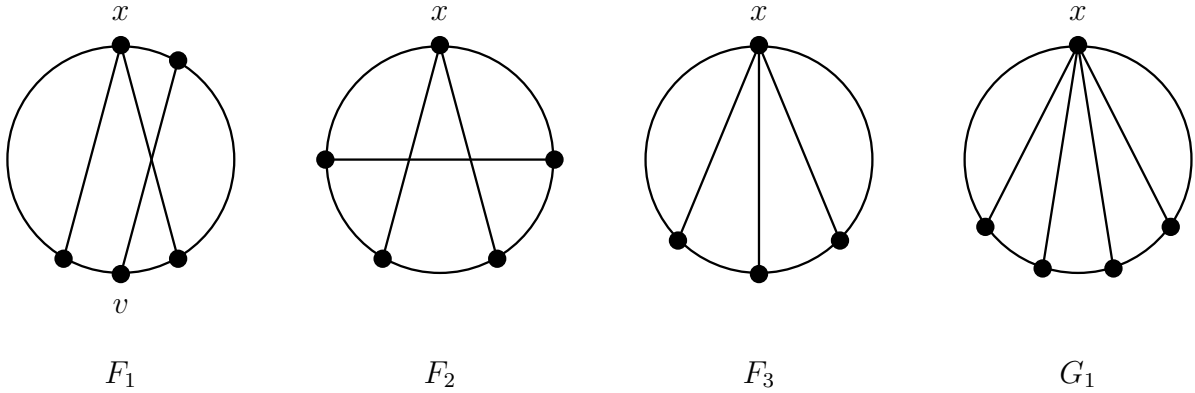

Next we consider the children of $E_2$.  The only multigraphs not
ruled out by Lemma~\ref{pairs} are the first three listed in
Figure~\ref{E2d}.  As
before, we need only consider the case where there is a cutvertex $x$ at
the vertex as indicated.  It can be checked that
$F_1$ can't lead to an
extremal graph by using the vertex marked $v$ in Lemma~\ref{con}.
For $F_2$, first note that if any of the edges containing $x$ is subdivided, and it doesn't
matter which one by symmetry, then taking $v$ to be the other
endpoint of that edge (after subdivision) in Lemma~\ref{con} shows that
$G$ is not extremal.  If none of the edges containing $x$ are
subdivided then $b:=|V(B)|$ satisfies $5\le b\le9$ because $G$
has Property $(\Delta)$, and applying
the $m$-bound to the vertex $x$ gives
$$
m(G)\le\max_{5\le b\le9}\{c(b-1,1)g(n-b)+g(n-b-1)\}<m_0.
$$
Finally, $F_3$ is treated the same way as $F_2$, 
noting that the two pairs of doubled edges must both be
subdivided in the same manner, the only edge containing $x$ which can
be subdivided further is the vertical one in the diagram, and the
maximum is now taken over $6\le b\le 9$.

The only grandchild of $E_2$ not thrown out by either Lemma~\ref{312} or
Lemma~\ref{pairs} is the multigraph $G_1$ in Figure~\ref{E2d}, a
child of $F_3$.  It is handled in the same way as $F_2$ and $F_3$ and
the reader should be able to fill in the details at this point.
It is easy to check that the children of $G_1$ are all eliminated, and
so we have finished with the descendants of $E_2$.

Lemma~\ref{pairs} rules out subdivisions of $E_3$ directly as well as,
in conjunction with Lemma~\ref{312}, many of its children and all of its
grandchildren.  The only surviving multigraphs not previously considered
are
those children listed in Figure~\ref{E3c}.  In $F_4$, we are reduced in
the usual manner to
the case where the vertex of degree at least four is a cutvertex $x$.
But then we can take $v$ as indicated in Lemma~\ref{con} and so this
child is not extremal.  In $F_5$ we need only consider when
there is a cutvertex $x$ in the doubled edge. But then
either $v_1$ or $v_2$ (depending on the placement of $x$) can be used in
Lemma~\ref{con} to take care of this child.  Similarly, in $F_6$ it is
easy to see by symmetry that no matter where the cutvertex is placed,
there is a $v$ for Lemma~\ref{con}.

\begin{figure}
\begin{center}
\begin{tabular}{ccccc}
\psset{xunit=0.012in, yunit=0.012in}
\psset{linewidth=1.0\psxunit}
\begin{pspicture}(-50,-64)(50,64)
\pscircle(0,0){50\psxunit}      
\psline(25, -43.3013)(25, 43.3013)      
\psline(-25, -43.3013)(-25, 43.3013)    
\psline(25, 0)(-25, 43.3013)    
\pscircle*(25, 43.3013){4\psxunit}      
\pscircle*(-25, 43.3013){4\psxunit}     
\pscircle*(25, -43.3013){4\psxunit}     
\pscircle*(-25, -43.3013){4\psxunit}    
\pscircle*(25, 0){4\psxunit}    
\rput[c](-25,57.3013){$x$}
\rput[c](25,-57.3013){$v$}
\end{pspicture}
&\rule{10pt}{0pt}&
\psset{xunit=0.012in, yunit=0.012in}
\psset{linewidth=1.0\psxunit}
\begin{pspicture}(-50,-64)(50,64)
\pscircle(0,0){50\psxunit}      
\psline(25, -43.3013)(25, 43.3013)      
\psline(-25, -43.3013)(-25, 43.3013)    
\psline(50, 0)(25, 0)   
\pscircle*(25, 43.3013){4\psxunit}      
\pscircle*(-25, 43.3013){4\psxunit}     
\pscircle*(25, -43.3013){4\psxunit}     
\pscircle*(-25, -43.3013){4\psxunit}    
\pscircle*(25, 0){4\psxunit}    
\pscircle*(50, 0){4\psxunit}    
\rput[c](25,57.3013){$v_1$}
\rput[c](25,-57.3013){$v_2$}
\end{pspicture}
&\rule{10pt}{0pt}&
\psset{xunit=0.012in, yunit=0.012in}
\psset{linewidth=1.0\psxunit}
\begin{pspicture}(-50,-64)(50,64)
\pscircle(0,0){50\psxunit}      
\psline(25, -43.3013)(25, 43.3013)      
\psline(-25, -43.3013)(-25, 43.3013)    
\psline(-25, 0)(25, 0)  
\pscircle*(25, 43.3013){4\psxunit}      
\pscircle*(-25, 43.3013){4\psxunit}     
\pscircle*(25, -43.3013){4\psxunit}     
\pscircle*(-25, -43.3013){4\psxunit}    
\pscircle*(25, 0){4\psxunit}    
\pscircle*(-25, 0){4\psxunit}   
\end{pspicture}
\\\\
$F_4$&&$F_5$&&$F_6$
\end{tabular}
\caption{Three children of $E_3$}\label{E3c}
\end{center}
\end{figure}
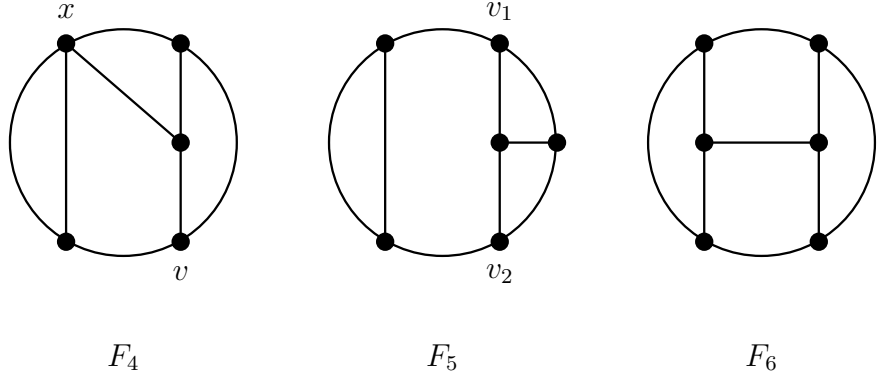

Finally we come to $E_4\iso K_4$.  Lemma~\ref{con} shows that if
$B$ is a subdivision of $E_4$ then we
need only consider when one of the degree 3 vertices is a cutvertex $x$.
If one of the edges  $xv$ is subdivided then $v$
satisfies Lemma~\ref{con}.  Therefore if there is a degree 2 vertex
in $B$, it must be formed by subdividing an edge between two
non-cutvertices of $E_4$, say $u$ and $v$.  Hence $G-\{u,v\}$ contains only one
nontrivial connected component with $n-3$ vertices, and we get
\begin{eqnarray*}
m(G)
&\le&
m(G-u)+m(G-N[u])\\
&\le&
m(G-\{u,v\})+m(G-u-N[v])+m(G-N[u])\\
&\le&
c(n-3)+g(n-5)+g(n-4)\\
&<&
m_0,
\end{eqnarray*}
so these graphs are not extremal.

We are reduced to considering the case
where $B\iso K_4$, that is, when no edges are subdivided.  By our
choice of $B$ and the cases we have disposed
of so far, we can assume that all other endblocks containing the
cutvertex $x$ are isomorphic to $K_2$, $K_3$, or $K_4$
Assume that there are $i$ copies of $K_3$ and $j$ copies of $K_4$
other than $B$.  If these are the only blocks of $G$, then
$m(G)=2^i3^{j+1}+1$ where $2i+3j+4=n=3r_0+2$.  This quantity is maximized
when $i=2$ and $j=r_0-2$, giving $m(G)=4\cdot 3^{r_0-1}+1<m_0$.

Hence we may assume that $G$ has other blocks.  We subdivide this case into three subcases.  First, if $j\ge 2$ (in other words, if $x$ lies in at least three
$K_4$ endblocks) then applying $m$-bound shows that
$$
m(G)\le 27g(n-10,r-21) + g(n-11,r-21)<m_0.
$$
Now consider $j=1$.  Here our upper bound is
\begin{eqnarray*}
m(G)
&\le&
9\cdot 2^i c(n-2i-7,r-i-14)+g(n-2i-8)\\
&\le&
9\cdot 2^i c(3r_0-2i-5,r_0-i-3)+g(3r_0-2i-6).
\end{eqnarray*}
This is a decreasing function of $i$ within each congruence class modulo $2$, and using this fact it is routine to check that $m(G)<m_0$.

We are left with the case where $j=0$.  The $m$-bound gives
\begin{eqnarray*}
m(G)
&\le&
3\cdot 2^i c(n-2i-4, r-i-7) + g(n-2i-5)\\
&\le&
3\cdot 2^ic(3r_0-2i-2,r_0-i+4)+g(3r_0-2i-3).
\end{eqnarray*}
Note that $3r_0-2i-2\ge 3(r_0-i+4)$ only for $i\ge 14$.  For these values of $i$, the upper bound is a decreasing function of $i$ within each congruence class modulo $2$,  so we only need to verify that $m(G)<m_0$ for $i\le 15$.  These cases can all be routinely checked, although when $i\Cong 1,2\ (\Mod 3)$ the desired inequality will hold only for sufficiently large $r_0$, and one must note that these cases can arise only for such sufficiently large $r_0$.  In particular, $n\ge 2i+5$ by our assumptions, and this implies that $r_0\ge 2i/3+1$.

The only child of $E_4$ that is not a child of any other $E_k$
and is not ruled out by
Lemma~\ref{312} is $F_7$, shown in Figure~\ref{E4c}.
One can verify by considering
several cases that whether or not there is a cutvertex in $B$,
there is a vertex $v$ satisfying the hypotheses of Lemma~\ref{con}.
Finally, the grandchildren of $E_4$ all fall under the purview of
Lemma~\ref{312}.

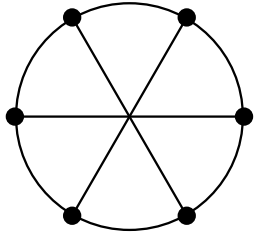
\begin{figure}
\begin{center}
\begin{center}
\psset{xunit=0.012in, yunit=0.012in}
\psset{linewidth=1.0\psxunit}
\begin{pspicture}(-50,-64)(50,64)
\pscircle(0,0){50\psxunit}      
\psline(25, -43.3013)(-25, 43.3013)     
\psline(-25, -43.3013)(25, 43.3013)     
\psline(50, 0)(-50, 0)  
\pscircle*(25, 43.3013){4\psxunit}      
\pscircle*(-25, 43.3013){4\psxunit}     
\pscircle*(25, -43.3013){4\psxunit}     
\pscircle*(-25, -43.3013){4\psxunit}    
\pscircle*(-50, 0){4\psxunit}   
\pscircle*(50, 0){4\psxunit}    
\end{pspicture}
\end{center}
\caption{$F_7$, a child of $E_4$}\label{E4c}
\end{center}
\end{figure}

We have now considered all the cases and so completed the proof of the
theorem.
\end{proof}

\section{Maximum independent sets}\label{sec-maximum}

We now turn to the consideration of maximum independent sets.  An
independent set $I$ if $G$ is {\it maximum\/} if it has maximum
cardinality over all independent sets of $G$.  We let
$m'(G)$ denote the number of maximum independent sets of $G$.  Since
every maximum independent set is also maximal we have
$m'(G)\le m(G)$, so for any finite family of graphs,
\begin{equation}
\label{mm'}
\max_{G\in\cF} m'(G)\le \max_{G\in\cF} m(G).
\end{equation}

We say that $G$ is {\it well covered\/} if every one of its maximal
independent sets is also
maximum.  Then we have equality in~(\ref{mm'}) if and only if some
graph with a maximum number of maximal independent sets is well 
covered.  The graphs $G(n)$, $G'(n)$, $C(n)$, and $G(n,r)$ are well covered for
all pairs $n,r$ and $C(n,r)$ is well covered when $n\Cong r\ (\Mod 2)$
and $n\ge 3r$, so we immediately have the following result.
\begin{theorem}\label{maximum-first-theorem}
Let $G$ be a graph with $n$ vertices and at most $r\ge 1$ cycles.
\begin{enumerate}
\item[(I)]  For all such graphs,
$$
m'(G)\le g(n,r),
$$
with equality if and only if $G\iso G(n,r)$.
\item[(II)]  If $n\ge3r$ and $n\Cong r\ (\Mod
2)$, or if $r\ge r(C(n))$, then for all such graphs that are connected,
$$
m'(G)\le c(n,r),
$$
with equality if and only if $G\iso C(n,r)$.
\end{enumerate}
\end{theorem}

This leaves only the case where $G$ is connected and $n\not\Cong r\ (\Mod 2)$.  To state and motivate the result in this case, we recall the work done on trees and graphs with at most one cycle.

When $n$ is even, the tree $K_1 * (K_1 \uplus \frac{n-2}{2}K_2)$ has $t(n)$ maximum independent sets (in fact, by the upcoming Theorem~\ref{zitotree}, it is the only such tree).  When $n$ is odd there is only one extremal tree for the maximal independent set problem, and it is not well-covered.  Define a family of trees by
$$
T'(n):=\case{K_1 * (K_1 \uplus \frac{n-2}{2}K_2)}
{if $n\Cong 0\ (\Mod 2)$,}
{K_1 * (2K_1 \uplus \frac{n-3}{2}K_2)}
{if $n\Cong 1\ (\Mod 2)$,\rule{0pt}{20pt}}
$$
and let
$$
t'(n):=m'(T'(n))=\case{2^\frac{n-2}{2}+1}
{if $n\Cong 0\ (\Mod 2)$,}
{2^\frac{n-3}{2}}
{if $n\Cong 1\ (\Mod 2)$.\rule{0pt}{20pt}}
$$
Zito~\cite{zit:smn} proved the following result.

\begin{theorem}[Zito~\cite{zit:smn}]
\label{zitotree}
If $T$ is a tree on $n\ge 2$ vertices  then
$$
m'(T)\le t'(n)
$$
with equality if and only if $T\iso T'(n)$.
\end{theorem}

In \cite{jc:nmi}, Jou and Chang gave a short proof of this theorem 
and considered graphs with at most one cycle.

\begin{theorem}[Jou and Chang~\cite{jc:nmi}]
\label{jc}
Let $G$ be a graph with at most one cycle and $n$ vertices where
$n\Cong 0\ (\Mod 2)$.  Then
$$
m'(G)\le t(n)
$$
with equality if and only if $G$ is a tree, and thus, by Theorem~\ref{zitotree}, if and only if $G\iso T'(n)$.
\end{theorem}

The rest of this section will be devoted to the proof of the maximum independent set version of Theorem~\ref{main1}.
\begin{theorem}
\label{main2}
Let $G$ be a connected graph with $n\ge 3r$ vertices and at most $r$
cycles, where $r\ge 1$ and $n\not\Cong r\ (\Mod 2)$.  Then
$$
m'(G)\le c(n,r-1),
$$
with equality if and only if $G$ has precisely $r-1$ cycles, and thus, by Theorem~\ref{main1}, if and only if $G\cong C(n,r-1)$ or if $G$ is isomorphic to one of the exceptional graphs listed there.
\end{theorem}

Before proving Theorem~\ref{main2} we present two lemmas.  The
first is an analogue the $m$-bound and $m$-recursion that proved useful
for maximal independent sets.  Its proof is similar and so is omitted.
\begin{lemma}
The invariant $m'(G)$ satisfies the following inequalities.
\begin{enumerate}
\item[(1)]\
If $v\in V$ then
$$
m'(G)\le m'(G-v)+m'(G-N[v]).
$$
\item[(2)]
If $G$ has a complete subgraph $B$ with at least one vertex
adjacent only to other vertices of $B$ then
$$
m'(G)\le\sum_{v\in V(B)} m'(G-N[v]).
$$
\end{enumerate}
\end{lemma}
We will refer to parts (1) and (2) of this lemma as the $m'$-bound and $m'$-recursion, respectively.

Our second lemma will be useful for eliminating the cases where a vertex 
serves as a cutvertex for more than one endblock.

\begin{lemma}\label{leaves}
Suppose that the graph $G$ contains a vertex $x$ and a set of at least
two other vertices $U$ such that the induced graph $G[U]$ is not
complete and for all $u\in U$,
\begin{equation}
\label{xU}
x\in N(u)\subseteq U\cup\{x\}.
\end{equation}
Then $m'(G)=m'(G-x)$.
\end{lemma}
\begin{proof}
It suffices to show that no maximum independent set of $G$ contains
$x$.  Suppose not, and let $I$ be a maximum independent set with $x\in
I$.  Then by~(\ref{xU}) we have $I\cap U=\emptyset$.  Also, since $G[U]$ is
not complete there is an independent set $A\subset U$ containing at least
two vertices.  But then $I\cup A-x$ is a larger independent set than
$I$, a contradiction.
\end{proof}

\medskip

\newenvironment{proof-main2}{\noindent {\bf Proof (of Theorem~\ref{main2}):}}{{\qed}}
\begin{proof-main2}
We will use induction on $r$.  The base case of $r=1$ is precisely
Theorem~\ref{jc}, so we will assume $r\ge 2$.  Note that since $n\ge 3r$ and
$n\not\Cong r\ (\Mod 2)$ we have $n\ge 3r+1$.  Let $G$ be an extremal graph
satisfying the hypothesis of the theorem.

If $G$ has less than $r$ cycles then we are done by
Theorem~\ref{main1}, so we may assume that $G$ has exactly $r$
cycles.  Since we have chosen $G$ to be extremal, we may assume that
$$
m'(G)\ge c(n,r-1)=3^{r-2}\cdot 2^{\frac{n-3r+5}{2}} + 2^{r-2}.
$$

Let $B$ be an endblock of $G$.  First, if $B$ has intersecting cycles,
then the Ear Decomposition Theorem shows that
$B$ must contain a subdivision of the multigraph $D$ (shown in Figure~\ref{D}).
This implies that $B$ contains a non-cutvertex of degree at least $3$
that lies in at least $3$ cycles, from which
$m'$-bound gives the contradiction
$$
m'(G)\le c(n-1,r-3)+g(n-4,r-3)<c(n,r-1).
$$
Hence $B$ is either $K_2$, $K_3$, or $C_p$ for some
$p\ge 4$.  Since these possibilities have at most one cycle
and we are assuming that $G$ has $r\ge 2$ cycles, $G$ cannot be a single
block.  Hence $B$ must contain a cutvertex $x$ of $G$.

First suppose that $B\iso C_p$ for some $p\ge 4$.
Label the vertices of
$B$ as $x,u,v,w,\ldots$ so that they read one of the possible directions
along the cycle.  Since $G-v$ is connected with $n-1$ vertices and $r-1$ cycles, 
induction applies to
give $m'(G-v)\le c(n-1,r-2)$.  Furthermore, $G-v$ has exactly $r-1$ cycles, so
by induction we cannot have equality.
Similarly, $G-N[v]$ has $n-3$ vertices and $r-1$ cycles, so $m'(G-N[v])<c(n-3,r-2)$.
An application of the $m'$-bound gives the contradiction
$$
m'(G)\le m'(G-v)+m'(G-N[v])<c(n-1,r-2)+c(n-3,r-2)=c(n,r-1).
$$

We now know that all endblocks of $G$ must be copies of either $K_2$ or 
$K_3$.   We claim that such endblocks must be disjoint.
Suppose to the  contrary that two endblocks share a vertex, which must
therefore be the cutvertex  $x$.  Considering the two cases when at 
least
one endblock is a $K_2$ (so that $G-x$ has an isolated vertex which
must be in each of its maximum independent sets) or when
both  are copies of $K_3$, we can use  Lemma~\ref{leaves} and the fact that $n\ge 3r+1$ to get
\begin{eqnarray*}
m'(G)	&=& m'(G-x)\\
	&\le&
\case{g(n-2,r)}{if one endblock is a $K_2$,}
{g(n-1,r-2)}{if both endblocks are $K_3$'s}
\\
	&<& c(n,r-1).
\end{eqnarray*}
This contradiction proves the claim.

Now let $B$ be an endblock of $G$ satisfying
Proposition~\ref{endB}, so that it intersects at most one non-endblock.
By what we have just shown, $B$ intersects precisely one
other block and that block is not an endblock.  We claim that
this block is isomorphic to $K_2$.  Suppose not.  Then the cutvertex
$x$ of $B$ lies in at least one cycle not contained in $B$ and is
adjacent to at least two vertices not in $B$.  Again, we consider the
cases $B\iso K_2$ and $B\iso K_3$ separately to obtain
\begin{eqnarray*}
m'(G)	&\le& m'(G-x)+m'(G-N[x])\\
	&\le&
\case{c(n-2,r-1)+g(n-4,r-1)}{if $B\iso K_2$,}
{2c(n-3,r-2)+g(n-5,r-2)}{if $B\iso K_3$}\\
	&<&c(n,r-1),
\end{eqnarray*}
proving our claim.

Since $G$ is not itself a block, Proposition~\ref{endB} shows that
$G$ contains at least two endblocks, say $B$ and $B'$, that each 
intersect at most one non-endblock.  Let the cutvertices of these 
endblocks be labeled $x$ and $x'$, respectively.  We have
shown that $B$ and $B'$ are disjoint and that they each intersect
precisely one other block, which must be isomorphic to $K_2$.
We claim that there is vertex $v_0$ such that these two
copies of $K_2$ have vertices $\{x,v_0\}$ and $\{x',v_0\}$.  If this is 
not the case then
\begin{equation}
\label{cup2}
|(N(x)-B)\cup (N(x')-B')|= 2.
\end{equation}

There are now three cases to consider depending on the nature of $B$
and $B'$.  First, suppose $B\iso B'\iso K_2$.  Using
the $m'$-recursion twice, we get
\begin{eqnarray*}
m'(G)	&\le&m'(G-B)+m'(G-N[x])\\
	&\le&m'(G-B)+m'(G-N[x]-B')+m'(G-N[x]-N[x']).
\end{eqnarray*}
Since the three graphs in this last expression may still have $r$ cycles, we will also have to use induction on $n$.  We consider two cases, depending on whether or not $G-N[x]-N[x']$ has parameters lying in the range of Theorem~\ref{main1}.  Clearly $G-B$ is connected with $n-2$ vertices and $r$ cycles, so if $n\ge 3r+5$ then $n-6\ge 3r-1$.  Also, $n-2$ and $r$ are of different parity with $n-2\ge 3r$.   
By assumption~(\ref{cup2}), $G-N[x]-N[x']$ has $n-6$ vertices and at most $r$ cycles yielding $m'(G-N[x]-N[x'])\le g(n-6,r)$.
Secondly, $G-N[x]-B'$ has $n-5$ vertices and at most $r$ cycles so $m'(G-N[x]-B')\le g(n-5,r)$.
Finally, we can apply induction to conclude that $m'(G-B)\le c(n-2,r-1)$.
Putting everything together we get
\begin{eqnarray*}
m'(G)	&\le&\case{c(n-2,r-1)+g(n-5,r)+g(n-6,r)}
			{if $n\ge3r+5$,}
			{c(n-2,r)+g(n-5)+g(n-6)}
			{if $n=3r+1$ or $n=3r+3$}\\
	&<&c(n,r-1),
\end{eqnarray*}
a contradiction.

Now suppose that $B\iso K_2$ and $B'\iso K_3$.  Proceeding in much the
same manner as before gives
\begin{eqnarray*}
m'(G)	&\le&m'(G-B)+2m'(G-N[x]-B')+m'(G-N[x]-N[x'])\\
	&\le&\case{c(n-2,r-1)+2g(n-6,r-1)+g(n-7,r-1)}
			{if $n\ge3r+3$,}
			{c(n-2,r)+2g(n-6)+g(n-7)}
			{if $n=3r+1$}\\
	&\le&c(n,r-1),
\end{eqnarray*}
and equality cannot occur because $G-B$ has exactly $r$ cycles.

The third case is when $B\iso K_3$ and $B'\iso K_3$.  Trying the same
technique we obtain
\begin{eqnarray*}
m'(G)	&\le&2m'(G-B)+2m'(G-N[x]-B')+m'(G-N[x]-N[x'])\\
	&\le&2c(n-3,r-2)+2g(n-7,r-2)+g(n-8,r-2)\\
	&\le&c(n,r-1).
\end{eqnarray*}
Again we cannot have equality throughout because $G-B$ has exactly $r-1$ cycles.

Now that we have established the existence of $v_0$, we are almost
done.  Observe that there is at most one block $C$ other than the
$K_2$'s connecting $v_0$ to endblocks and those endblocks themselves.
(If there
were more than one such block, then since endblocks can't intersect
this would force the existence of another $K_2$ and corresponding
endblock which we hadn't considered.)  So $C$, if it exists, must be
an endblock containing $v_0$.  By our characterization of endblocks,
this leaves only three possibilities, namely $C\iso \emptyset, K_2, K_3$.
It is easy to check the corresponding graphs $G$ either do not exist
because of parity considerations or satisfy
$m'(G)<c(n,r-1)$.  We have now shown that no graph with exactly $n$
vertices and $r$ cycles has as many maximum independent sets as
$C(n,r-1)$, and thus finished the proof of Theorem~\ref{main2}.
\end{proof-main2}

\section{The gap revisited}\label{sec-gap2}

We now need to look at  maximum independent sets in the gap.  Consider
first the case when $G$ is connected with $n=3r_0+2$ vertices
and less than $r_1$ cycles.
Then $n\Cong r_0\ (\Mod 2)$ and by Theorem~\ref{gap} we have
$m'(G)\le m(G)\le c(n,r_0)$, with the second inequality reducing to an
equality if and only if $G\iso C(n,r_0)$ or $C(n,r_0+1)$.  Since the
former graph is well covered but the latter is not, we have the
following result.

\begin{theorem}
\label{gap2}
Let $G$ be a connected graph with $n$ vertices, $n=3r_0+2$
where $n\ge 7$, and less than $r_1$ cycles.  Then
$$
m'(G)\le c(n,r_0)
$$
with equality if and only if $G\iso C(n,r_0)$.
\end{theorem}

For the $n=3r_0+1$ case we need to adapt the proof of Theorem~\ref{gap}.
The result mirrors the trend exhibited by Theorem~\ref{main2}.

\begin{theorem}
Let $G$ be a connected graph with $n$ vertices, $n=3r_0+1$
where $n\ge 7$, and less than $r_1$ cycles.  Then
$$
m'(G)\le c(n,r_0-1)
$$
with equality if and only if $G\iso C(n,r_0-1)$.
\end{theorem}
\begin{proof}  We will use induction on $n$.  The $n=7$ case
has been checked by
computer and so we assume $G$ is a graph satisfying the hypotheses of
the theorem where $n\ge 10$, or equivalently, $r_0\ge3$.  We
will begin as in the proof of Theorem~\ref{gap}, considering
possible endblocks produced by the inductive procedure in
the Ear Decomposition Theorem.

Our first order of business will be to show that any endblock of an
extremal $G$ must be isomorphic to $K_i$ for $2\le i\le4$ or the graph
$D_1$ shown in Figure~\ref{D_1}.  Note that, unlike in the proof of Theorem~\ref{gap},
here we will consider all endblocks of $G$, not just those that satisfy
the conclusion of Proposition~\ref{endB}.  As the proof of
Lemma~\ref{2deg2} no longer
holds for maximum independent sets, we need the following result to
replace Property ($\Delta$).  In it, and in the future, it will be
convenient to use the notation
$$
m_0'=c(n,r_0-1)=8\cdot 3^{r_0-2}+2^{r_0-2}.
$$
\begin{lemma}[Path Lemma]
\label{path}
Suppose there is a path $P=v_1v_2v_3$ in a block of $G$ satisfying the
following three conditions.
\begin{enumerate}
\item[(1)] $\deg v_1\ge3$, $\deg v_2=2$, and $v_1v_3\not\in E(G)$,
\item[(2)] $G-P$ is connected,
\item[(3)] One of the following two subconditions hold
	\begin{enumerate}
	\item[(a)]  $G-N[v_1]$ is connected, or
	\item[(b)]  $G-v_1-N[v_3]=G_1\uplus G_2$ where $G_1$ is
	connected and $|V(G_2)|\le2$.
	\end{enumerate}
\end{enumerate}
Then $m'(G)<m_0'$.
\end{lemma}
\begin{proof}  Using the $m'$-bound twice gives
\begin{equation}
\label{3m'}
m'(G)\le m'(G-v_1-v_3)+m'(G-v_1-N[v_3])+m'(G-N[v_1]).
\end{equation}
Conditions (1) and (2) of the current lemma imply that
$G-v_1-v_3=H\uplus\{v_2\}$ where $H$ is connected.  Furthermore, since
$P$ is in a block, $v_2$ must lie in at least one cycle of $G$.  Hence $H$
has $n-3$ vertices and less than $r_1-1$ cycles with these two parameters
satisfying the hypotheses of the theorem.  By induction,
$$
m'(G-v_1-v_3)\le c(n-3,r_0-2).
$$

Now suppose (3a) holds.  Then $m'(G-N[v_1])\le c(n-4)$.  Also $\deg
v_3\ge2$ (since it is in a non-$K_2$ block) and $v_1\not\in N[v_3]$ by
condition~(1), so $m'(G-v_1-N[v_3])\le g(n-4)$.  Putting all these bounds
into~(\ref{3m'}) gives
$$
m'(G)\le c(n-3,r_0-2)+c(n-4)+g(n-4)\le m_0'.
$$
Equality can only be achieved if $r_0=3$, $H\iso C(7,1)$, and
$G-v_1-N[v_3]\iso G(6)$.  But $C(7,1)$ has only one
cycle while $G(6)$ has two, contradicting the fact that
$G-v_1-N[v_3]\subset H$, so we must have $m'(G)<m_0'$ in this case.

Next we look at (3b).  Considering the cases where $|V(G_2)|=0$, 1,
or 2 gives
$$
m'(G-v_1-N[v_3])\le\max\{c(n-4),\ c(n-5),\ 2c(n-6)\}=c(n-4),
$$
while $m'(G-N[v_1])\le g(n-4)$.  Thus we get the same bound on $m'(G)$ as
in the case where (3a) held.  Equality implies $H\iso C(7,1)$ and
$G_1\iso C(6)$ or $K_4$ since $c(6)=2c(4)$, but then we have the
same problem with cycles.  This final contradiction ends the proof of
Lemma~\ref{path}.
\end{proof}

In all of our applications of the Path Lemma we will set up the
notation so that $v_1=v$.

\begin{lemma}
\label{pathco}
Let $B$ be a block of $G$ which comes from subdividing a multigraph $M$, and suppose
that $M$ and $B$ satisfy either
\begin{enumerate}
\item[(1)] $M$ contains an edge $uv$ such that
\begin{enumerate}
\item[(a)] $v$ is a non-cutvertex in $G$,
\item[(b)] $\deg_M v\ (=\deg_B v)=3$, and
\item[(c)] $uv$ is subdivided more than twice in $B$ and none of these inserted vertices are cutvertices,
\end{enumerate}
or
\item[(2)] $M$ contains a doubled edge where both edges are subdivided exactly once and neither of these inserted vertices is a cutvertex.
\end{enumerate}
Then $G$ is not extremal.
\end{lemma}
\begin{proof}
For (1), suppose such an  edge $vw$ is subdivided
three or more times.  Then the hypotheses  of the Path
Lemma are satisfied with
condition (3b) and $v_1=v$, so $m'(G)<m_0'$.

Part (2) of the lemma follows from the fact that extremal graphs
cannot have duplicated vertices: if both edges are subdivided
exactly once then we have a pair
$s,t$ of duplicated vertices, and Lemma~\ref{dup} shows that
$$
m'(G)\le m(G)= m(G-s)\le c(n-1)<m_0',
$$
another contradiction.
\end{proof}

We are now ready to begin restricting the type of endblocks an extremal
graph may possess.  Let $B$ denote an endblock of $G$.
First we consider the case where we have
$B\iso C_p$.  As $c(n,1)<m_0'$
for all $n\ge 10$, $G$ may not itself be a block, and thus we may
assume that there is a cutvertex, say $x$, in $B$.
The $p=4$ case is ruled out by Lemma~\ref{pathco}~(2).  If $p=5$ then
considering the cases where $x$ is adjacent to exactly one or
more than one vertex outside $B$ gives
\begin{eqnarray*}
m'(G)	&\le& m'(G-x)+m'(G-N[x])\\
	&\le&
\case{3c(n-5)+2g(n-6)}{if $x$ has exactly one neighbor outside $B$,}
{3g(n-5)+2g(n-7)}{if $x$ has more than one neighbor outside $B$}\\
	&<&m_0'
\end{eqnarray*}
If $p\ge 6$, then label the non-cutvertices of $B$
by $v_1,v_2,\dots$ so that $xv_1v_2\dots v_{p-1}$ is a path.
Using the $m'$-bound twice gives
\begin{eqnarray*}
m'(G)
&\le&
m'(G-v_2)+m'(G-N[v_2])\\
&\le&
m'(G-v_2-v_4)+m'(G-v_2-N[v_4])+m'(G-N[v_2]).
\end{eqnarray*}
We can apply induction to $m'(G-v_2-v_4)$ since it consists of an isolated
vertex together with a connected graph with $n-3$ vertices and one less cycle
than $G$, and similarly to $m'(G-N[v_2])$.  This gives
$$
m'(G)\le 2c(n-3,r_0-2)+c(n-4)\le m_0',
$$
with equality if and only if $r_0=3$.  But then $G-v_2-N[v_4]\iso C(6)$,
which contradicts the fact that the former graph has a vertex of
degree 1 while the latter does not.
So if $B$ is a cycle in an extremal graph then $B\iso K_3$.

We now consider the case where $B$ comes from subdividing the graph
$D$ in Figure~\ref{D}.  If $G$ is itself
a subdivision of $D$ then Lemma~\ref{pathco}~(1) shows that $G$ has
at most 8 vertices, which falls within the range of the computer calculations
we have performed.  Therefore $B$ must contain a cutvertex $x$.

First
suppose that $x$ is a vertex of $D$ and let $v$ be the other vertex of
degree 3.  If all of the edges of $D$ are subdivided
then at most one of them can be subdivided once by
Lemma~\ref{pathco}~(2), and so $G$ is not
extremal because of option (3a) in the Path Lemma.  If
one of the edges is not subdivided, then the only two possibilities
for the number of subdividing vertices in the three edges are $(2,1,0)$
and $(2,2,0)$ by Lemma~\ref{pathco}~(1).  It is then
easy to check that
\begin{eqnarray*}
m'(G)
&\le&
m'(G-x)+m'(G-N[x])\\
&\le&
\case
{3g(n-5)+g(n-6)}{in the $(2,1,0)$ case}
{g(n-6)+g(n-7)}{in the $(2,2,0)$ case,}
\\
&<&m_0'.
\end{eqnarray*}

This leaves the case where $x$ is interior to an edge
of $D$.  Neither of the other edges of $D$ may be subdivided more than
twice by Lemma~\ref{pathco}~(1), so if one of them is subdivided twice then
$G$ is not extremal by the Path Lemma~(3b).
Lemma~\ref{pathco}~(2) rules out
the case where both edges not containing $x$ are subdivided once.
Therefore we may assume that one of these edges is subdivided once
and the other is not subdivided.  Furthermore, the Path Lemma~(3a) 
shows that the edge containing $x$ may be subdivided at most once
to either side of $x$.  Thus we are reduced to considering the three
endblocks $D_1$, $D_2$, and $D_3$ shown in Figure~\ref{D_1}.

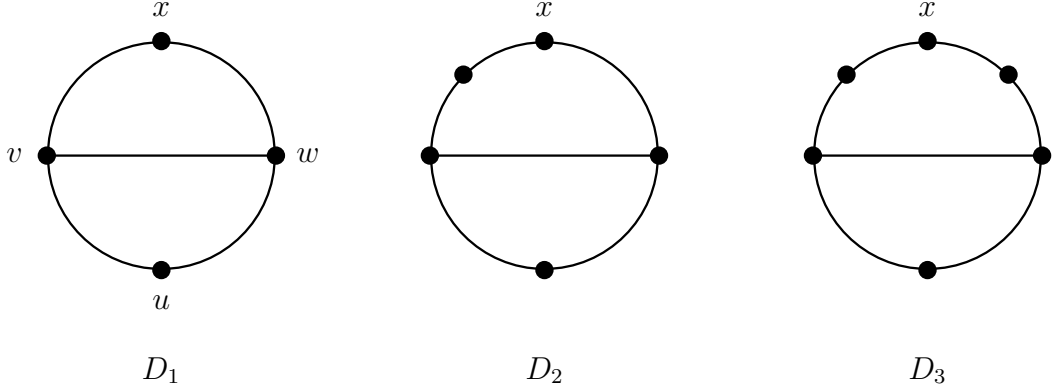
\begin{figure}
\begin{center}
\begin{tabular}{ccccc}
\psset{xunit=0.012in, yunit=0.012in}
\psset{linewidth=1.0\psxunit}
\begin{pspicture}(-64,-64)(64,64)
\pscircle(0,0){50\psxunit}      
\psline(50, 0)(-50, 0)  
\pscircle*(0, 50){4\psxunit}    
\pscircle*(0, -50){4\psxunit}   
\pscircle*(-50, 0){4\psxunit}   
\pscircle*(50, 0){4\psxunit}    
\rput[c](-64,0){$v$}
\rput[c](0,64){$x$}
\rput[c](64,0){$w$}
\rput[c](0,-64){$u$}
\end{pspicture}
&\rule{10pt}{0pt}&
\psset{xunit=0.012in, yunit=0.012in}
\psset{linewidth=1.0\psxunit}
\begin{pspicture}(-64,-64)(64,64)
\pscircle(0,0){50\psxunit}      
\psline(50, 0)(-50, 0)  
\pscircle*(0, 50){4\psxunit}    
\pscircle*(0, -50){4\psxunit}   
\pscircle*(-50, 0){4\psxunit}   
\pscircle*(50, 0){4\psxunit}    
\pscircle*(-35.3553, 35.3553){4\psxunit}        
\rput[c](0,64){$x$}
\end{pspicture}
&\rule{10pt}{0pt}&
\psset{xunit=0.012in, yunit=0.012in}
\psset{linewidth=1.0\psxunit}
\begin{pspicture}(-64,-64)(64,64)
\pscircle(0,0){50\psxunit}      
\psline(50, 0)(-50, 0)  
\pscircle*(0, 50){4\psxunit}    
\pscircle*(0, -50){4\psxunit}   
\pscircle*(-50, 0){4\psxunit}   
\pscircle*(50, 0){4\psxunit}    
\pscircle*(-35.3553, 35.3553){4\psxunit}        
\pscircle*(35.3553, 35.3553){4\psxunit} 
\rput[c](0,64){$x$}
\end{pspicture}
\\\\
$D_1$&&$D_2$&&$D_3$
\end{tabular}
\caption{Three special subdivisions of $D$}\label{D_1}
\end{center}
\end{figure}

Both $D_2$ and $D_3$ can be eliminated by using the
$m'$-bound on the vertex $x$:
\begin{eqnarray*}
m'(G)
&\le&
m'(G-x)+m'(G-N[x])\\
&\le&
\case{2g(n-5)+2g(n-6)}{if $B\iso D_2$,}
{g(n-6)+3g(n-7)}{if $B\iso D_3$}\\
&<&
m_0'.
\end{eqnarray*}
This leaves us with the case $B\iso D_1$, which we postpone until later.

We now need to go through the $E$ graphs from Figure~\ref{E}.  It will
be useful to have an analogue of the large degree and large number of
cycles results (Lemmas~\ref{36} and \ref{312}) used in proving the first Gap Theorem.
\begin{lemma}\label{bound}
Suppose that $v\in V(G)$ is a non-cutvertex satisfying either
\begin{enumerate}
\item[(1)] $\deg v\ge 3$ and $r(v)\ge 7$, or
\item[(2)] $\deg v\ge 4$ and $r(v)\ge 6$.
\end{enumerate}
Then $m'(G)<m_0'$.
\end{lemma}
\begin{proof}
Since $r\le r_0+5$, we have
\begin{eqnarray*}
m'(G)	&\le&m'(G-v)+m'(G-N[v]),\\
	&\le&
\case{c(n-1,r_0-2)+g(n-4,r_0-2)}{if $\deg v\ge 3$ and $r(v)\ge 7$,}
{c(n-1,r_0-1)+g(n-5)}{if $\deg v\ge 4$ and $r(v)\ge 6$}\\
	&<&m_0',
\end{eqnarray*}
proving the lemma.
\end{proof}

No subdivision of $E_1$ can be an endblock in an extremal graph 
by part (2) of the
previous lemma.
The same reasoning shows that $B$ may only be a subdivision of $E_2$
if the vertex of degree 4 in $E_2$ is the cutvertex $x$.
Label the other two
vertices $v$ and $w$.  Lemma~\ref{pathco} shows that the number of
subdividing vertices for the pair of $vx$ edges must be one of
$(1,0)$, $(2,0)$, $(2,1)$, or $(2,2)$.  The last two cases are
eliminated by~(3a) of the Path Lemma.  The same
can be said of the $wx$ edges and that $vw$ can be subdivided at most
once by the Path Lemma~(3b).  So $5\le|V(B)|\le8$
and
\begin{eqnarray*}
m'(G)	&\le&m'(G-x)+m'(G-N[x])\\
	&\le&\max_{5\le b\le 8}\{t'(b-1)g(n-b)\}+g(n-6)\\
	&<&m_0',
\end{eqnarray*}
where we remind the reader that $t'(n)$ denotes the maximum number of maximum independent sets in a tree with $n$ vertices.

To deal with $E_3$, we need an analogue of the Triangle Lemma in this
setting.
\begin{lemma}[Strict Triangle Lemma]
Suppose $G$ contains three vertices $\{v_1,v_2,v_3\}$ satisfying the
following two restrictions.
\begin{enumerate}
\item[(1)]  These vertices form a $K_3$ with $\deg v_2=2$ and
	$\deg v_1, \deg v_3\ge3$.
\item[(2)] $G-N[v_2]$ is connected and at least one of $G-N[v_1]$ or $G-N[v_3]$
	is connected.
\end{enumerate}
Then $m'(G)<m_0'$.
\end{lemma}
\begin{proof}
Since $G-N[v_2]$ has $n-3$ vertices and less than $r_1-1$ cycles, we can
use induction to conclude $m'(G-N[v_2])\le c(n-3,r_0-2)$.  Now using
the $m'$-recursion
\begin{eqnarray*}
m'(G)	&\le&m'(G-N[v_2])+m'(G-N[v_1])+m'(G-N[v_3])\\
	&\le&c(n-3,r_0-2)+c(n-4)+g(n-4)\\
	&\le&m_0'.
\end{eqnarray*}
Equality forces $r_0=3$, $G-N[v_2]\iso C(7,1)$, and $G-N[v_1]\iso
C(6)$ or  $G(6)$.  Considering numbers of cycles and containments gives a
contradiction.
\end{proof}

Now consider subdivisions of $E_3$.  Let $vw$ be one of the doubled
edges.  By symmetry, we can assume that
the cutvertex of $B$ in $G$ (if there is one) is neither
in one of the $vw$ edges nor adjacent to $v$.  If one of the $vw$
edges is subdivided more than
once then $G$ is not extremal by (3b) of the Path Lemma and
Lemma~\ref{pathco}~(1).  By
Lemma~\ref{pathco}~(2) the only other option is to have one edge
subdivided once and the other not subdivided at all. But then the
Strict Triangle Lemma shows that $G$ is not extremal.

Finally we come to $E_4\iso K_4$.  First we claim that any edge
$vw$ of $K_4$ that does not contain a cutvertex of $G$ cannot
be subdivided.  By
the Path Lemma~(3b) such an edge cannot be subdivided more than once.
If it is subdivided exactly once, then we can use the $m'$-bound
twice and induction to get
\begin{eqnarray*}
m'(G)	&\le&m'(G-v-w)+m'(G-v-N[w])+m'(G-N[v])\\
	&\le&c(n-3,r_0-2)+g(n-5)+g(n-4)\\
	&<&m_0'.
\end{eqnarray*}
So in order for $G$ to have at least 10 vertices this $K_4$ must contain a cutvertex $x$ of $G$.
Suppose first that $x$ is a vertex of $K_4$ (before subdivision).
Let $u,v,w$ be the other three vertices of $K_4$.  Since none of
the edges between these three vertices are subdivided, we can use (3a)
of the Path Lemma to conclude that the edge $vx$ is subdivided at most
once.  If $vx$ is subdivided exactly once, then
\begin{eqnarray*}
m'(G)	&\le&m'(G-v)+m'(G-N[v])\\
	&\le&c(n-1,r_0-1)+c(n-4,r_0-2)\\
	&<&m_0'.
\end{eqnarray*}
If $x$ is interior to an edge of $K_4$, then taking $v$ to be a vertex of $K_4$ which is not adjacent to $x$ gives the same inequality, so $G$ is not extremal in this case either.  This shows that if $B$ comes from subdividing $K_4$ then we must have $B\iso K_4$ and one of the vertices of $B$ is a cutvertex.  We will return to eliminate this case at the end of the proof.

Like in the $n\Cong 1\ (\Mod 3)$ case of the first Gap Theorem,
one can use Lemma~\ref{bound} to rule out all descendants of the $E$ 
graphs.  So now we know that all endblocks
are copies of $K_i$, $2\le i\le 4$, or $D_1$.

The rest of our proof will parallel the last part of the demonstration of Theorem~\ref{main2}.  There we were able to show that any two endblocks are disjoint.  Here we will have to settle for showing that only copies of $D_1$ may intersect.

Suppse that the endlbocks $B$
and $B'$ both contain the cutvertex $x$.  If
$B\iso K_i$ and $B'\iso K_j$ where $2\le i,j \le4$ then
Lemma~\ref{leaves} shows that
$$
m'(G)=m'(G-x)\le (i-1)(j-1)g(n-i-j+1)<m_0'
$$
unless $i=j=4$.  But in that case removing $B$ and $B'$ destroys $14$ of the
at most $r_0+5$ cycles and we can use the bound $m'(G-x)\le 9g(n-7,r_0-9)<m_0'$ 
instead.

Next consider the case where $B\iso K_i$, $2\le i\le4$, and
$B'\iso D_1$ where $D_1$ is labeled as in Figure~\ref{D_1}.    Then using
the $m'$-bound, Lemma~\ref{leaves}, and induction,  we have
\begin{eqnarray*}
m'(G)	&\le& m'(G-u)+m'(G-N[u])\\
	&\le&
\left\{\begin{array}{ll}
2g(n-5)+c(n-3,r_0-2)&\mbox{if $B\cong K_2$,}\\
4g(n-6)+c(n-3,r_0-2)&\mbox{if $B\cong K_3$,}\\
6g(n-7,r_0-5)+c(n-3,r_0-2)&\mbox{if $B\cong K_4$,}
\end{array}
\right.\\
	&<&m_0'.
\end{eqnarray*}

Hence we know that if two or more endblocks intersect at a cutvertex
$x$ then they must all be copies of $D_1$.  Now choose $B$ so that it intersects at most one other
non-endblock (by Proposition~\ref{endB}).   We can eliminate
$K_4$ as a possibility for $B$ by
labeling the cutvertex $x$ and counting cycles:
\begin{eqnarray*}
m'(G)	&\le& m'(G-x)+m'(G-N[x])\\
	&\le& 3c(n-4,r_0-2)+g(n-5,r_0-2)\\
	&<& m_0'.
\end{eqnarray*}
If $B\cong D_1$ doesn't intersect another $D_1$ endblock then the Strict Triangle Lemma can be used to show that $G$ is not extremal.  If $B\cong D_1$ intersects with $i$ other endblocks isomorphic to $D_1$ then we can adapt the proof of the Strict Triangle Lemma to show that $G$ is not extremal.  Label the vertices of $B$ as in Figure~\ref{D_1}.  Using the $m'$-recursion we get
$$
m'(G)
\le
m'(G-N[u])+m'(G-N[v])+m'(G-N[w])\\
$$
If $i\ge 2$, then both $m'(G-N[v])$ and $m'(G-N[w])$ lie in the range of Theorem~\ref{main1}, and by applying induction to $G-N[u]$ we have
$$
m'(G)\le c(n-3,r_0-2)+2\cdot 3^ic(n-3i-4,r_0-3i+2).
$$
This is a decreasing function of $i$, and it is strictly less that $m_0'$ when $i=2$.

\begin{figure}
\begin{center}
\begin{tabular}{ccccc}
\psset{xunit=0.012in, yunit=0.012in}
\psset{linewidth=1.0\psxunit}
\begin{pspicture}(0,0)(127.011,94.5937)
\psline(101.931, 57.9287)(63.5056, 64.7041)     
\psline(88.5859, 94.5937)(63.5056, 64.7041)     
\psline(88.5859, 94.5937)(101.931, 57.9287)     
\psline(38.4253, 94.5937)(63.5056, 64.7041)     
\psline(25.0803, 57.9287)(63.5056, 64.7041)     
\psline(25.0803, 57.9287)(38.4253, 94.5937)     
\psline(127.011, 87.8183)(101.931, 57.9287)     
\psline(127.011, 87.8183)(88.5859, 94.5937)     
\psline(0, 87.8183)(38.4253, 94.5937)   
\psline(0, 87.8183)(25.0803, 57.9287)   
\psline(63.5056, 33.7906)(63.5056, 64.7041)     
\psline(43.9966, 0)(63.5056, 33.7906)   
\psline(83.0147, 0)(63.5056, 33.7906)   
\psline(83.0147, 0)(43.9966, 0) 
\pscircle*(63.5056, 64.7041){4\psxunit} 
\pscircle*(101.931, 57.9287){4\psxunit} 
\pscircle*(88.5859, 94.5937){4\psxunit} 
\pscircle*(38.4253, 94.5937){4\psxunit} 
\pscircle*(25.0803, 57.9287){4\psxunit} 
\pscircle*(127.011, 87.8183){4\psxunit} 
\pscircle*(0, 87.8183){4\psxunit}       
\pscircle*(63.5056, 33.7906){4\psxunit} 
\pscircle*(43.9966, 0){4\psxunit}       
\pscircle*(83.0147, 0){4\psxunit}       
\end{pspicture}
&\rule{10pt}{0pt}&
\psset{xunit=0.012in, yunit=0.012in}
\psset{linewidth=1.0\psxunit}
\begin{pspicture}(0,0)(117.054,101.372)
\psline(97.5452, 67.5813)(58.5271, 67.5813)     
\psline(78.0361, 101.372)(58.5271, 67.5813)     
\psline(78.0361, 101.372)(97.5452, 67.5813)     
\psline(39.0181, 101.372)(58.5271, 67.5813)     
\psline(19.509, 67.5813)(58.5271, 67.5813)      
\psline(19.509, 67.5813)(39.0181, 101.372)      
\psline(117.054, 101.372)(97.5452, 67.5813)     
\psline(117.054, 101.372)(78.0361, 101.372)     
\psline(0, 101.372)(39.0181, 101.372)   
\psline(0, 101.372)(19.509, 67.5813)    
\psline(39.0181, 33.7906)(58.5271, 67.5813)     
\psline(78.0361, 33.7906)(58.5271, 67.5813)     
\psline(78.0361, 33.7906)(39.0181, 33.7906)     
\psline(58.5271, 0)(39.0181, 33.7906)   
\pscircle*(58.5271, 67.5813){4\psxunit} 
\pscircle*(97.5452, 67.5813){4\psxunit} 
\pscircle*(78.0361, 101.372){4\psxunit} 
\pscircle*(39.0181, 101.372){4\psxunit} 
\pscircle*(19.509, 67.5813){4\psxunit}  
\pscircle*(117.054, 101.372){4\psxunit} 
\pscircle*(0, 101.372){4\psxunit}       
\pscircle*(39.0181, 33.7906){4\psxunit} 
\pscircle*(78.0361, 33.7906){4\psxunit} 
\pscircle*(58.5271, 0){4\psxunit}       
\end{pspicture}
&\rule{10pt}{0pt}&
\psset{xunit=0.012in, yunit=0.012in}
\psset{linewidth=1.0\psxunit}
\begin{pspicture}(0,0)(127.011,97.4708)
\psline(101.931, 60.8059)(63.5056, 67.5813)     
\psline(88.5859, 97.4708)(63.5056, 67.5813)     
\psline(88.5859, 97.4708)(101.931, 60.8059)     
\psline(38.4253, 97.4708)(63.5056, 67.5813)     
\psline(25.0803, 60.8059)(63.5056, 67.5813)     
\psline(25.0803, 60.8059)(38.4253, 97.4708)     
\psline(127.011, 90.6954)(101.931, 60.8059)     
\psline(127.011, 90.6954)(88.5859, 97.4708)     
\psline(0, 90.6954)(38.4253, 97.4708)   
\psline(0, 90.6954)(25.0803, 60.8059)   
\psline(63.5056, 45.0542)(63.5056, 67.5813)     
\psline(63.5056, 22.5271)(63.5056, 45.0542)     
\psline(63.5056, 0)(63.5056, 22.5271)   
\pscircle*(63.5056, 67.5813){4\psxunit} 
\pscircle*(101.931, 60.8059){4\psxunit} 
\pscircle*(88.5859, 97.4708){4\psxunit} 
\pscircle*(38.4253, 97.4708){4\psxunit} 
\pscircle*(25.0803, 60.8059){4\psxunit} 
\pscircle*(127.011, 90.6954){4\psxunit} 
\pscircle*(0, 90.6954){4\psxunit}       
\pscircle*(63.5056, 45.0542){4\psxunit} 
\pscircle*(63.5056, 22.5271){4\psxunit} 
\pscircle*(63.5056, 0){4\psxunit}       
\end{pspicture}
\\\\
$G_1$&&$G_2$&&$G_3$
\end{tabular}
\caption{Three graphs to eliminate}\label{Ex10}
\end{center}
\end{figure}
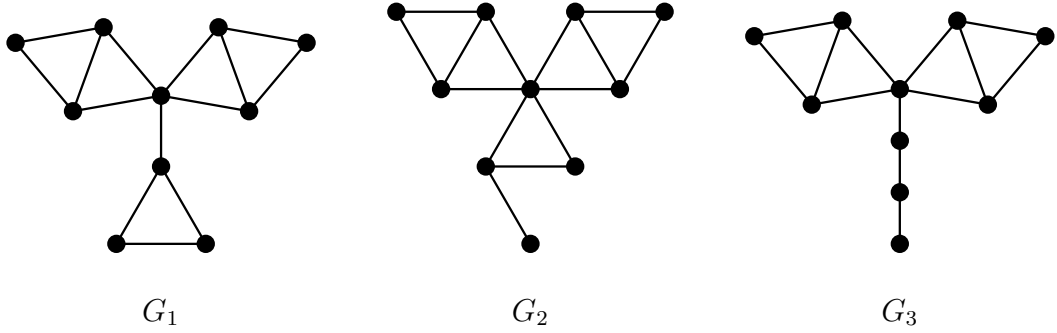

This leaves the $i=1$ case.  Here we use the bound
$$
m'(G-N[v])=m'(G-N[w])\le 3\cdot c(n-7)
$$
and induction to get
$$
m'(G)\le 20\cdot 3^{r_0-3}+7\cdot 2^{r_0-3}.
$$
This quantity is strictly less than $m_0'$ for $r_0\ge 4$.  As we are assuming $r_0\ge 3$, we have only to eliminate the case $r_0=3$.  In this case $G$ has $10$ vertices, of which $7$ are accounted for by the two intersecting copies of $D_1$.  This leaves $3$ vertices unaccounted for, and since we have limited the types of endblocks that can occur, there are only three possible graphs of this description.  These graphs are depicted in Figure~\ref{Ex10}.  It is easy to check that $m'(G_1)=2$, $m'(G_2)=10$, and $m'(G_3)=11$, which are all less than $c(10,2)=26$.

Hence if $B$ intersects precisely one non-endblock then $B$ is isomorphic to either $K_2$ or $K_3$, and thus by our previous work $B$ intersects precisely one other block, say $A$.
We claim that $A\iso K_2$ in the other two cases.  If not, then the
cutvertex $x$ is adjacent to two vertices of a cycle in $A$.
Using the $m'$-bound as well as induction
\begin{eqnarray*}
m'(G)	&\le& m'(G-x)+m'(G-N[x])\\
		&\le&
\case{c(n-2)+g(n-4)}{if $B\iso K_2$,}
{2c(n-3,r_0-2)+g(n-5)}{if $B\iso K_3$}\\
	&<& m_0'
\end{eqnarray*}
for $r_0\ge 3$, proving our claim that $A\iso K_2$.  
Thus if an endblock $B$ intersects at most one non-endblock $A$ then 
$B$ is
the only endblock intersecting $A$, $B\iso K_2$ or $K_3$, and $A\iso 
K_2$.

Suppose $B$ and $B'$ are endblocks of the type considered in the previous paragraph with cutvertices $x$ and $x'$, respectively.  We now claim that the
associated $K_2$ blocks must have vertex sets $\{x,v_0\}$ and 
$\{x',v_0\}$
for some $v_0$.  Suppose not and consider first the case $B\iso B'\iso 
K_2$.
The same argument as in Theorem~\ref{main2} shows that
\begin{eqnarray*}
m'(G)	&\le&m'(G-B)+m'(G-N[x]-B')+m'(G-N[x]-N[x'])\\
	&\le&c(n-2)+g(n-5)+g(n-6)\\
	&<&m_0'
\end{eqnarray*}
for $r_0\ge 3$, so such graphs are not extremal.
Now suppose that $B\iso K_3$ and $B'\iso K_2$.  In order to apply 
induction, it is important to use the $m'$-recursion
first on $B$ and then on $B'$ to get
\begin{eqnarray*}
m'(G)	&\le&2m'(G-B)+m'(G-N[x]-B')+m'(G-N[x]-N[x'])\\
	&\le&2c(n-3,r_0-2)+g(n-6)+g(n-7)\\
	&<&m_0',
\end{eqnarray*}
again resulting in a non-extremal graph.
Finally, if $B\iso B'\iso K_3$ then
\begin{eqnarray*}
m'(G)	&\le&2m'(G-B)+2m'(G-N[x]-B')+m'(G-N[x]-N[x'])\\
	&\le&2c(n-3,r_0-2)+2g(n-7)+g(n-8)\\
	&=&m_0'.
\end{eqnarray*}
Equality can only be achieved if $G-B\iso C(n-3,r_0-2)$.  But then $r(G-
B)=r_0-2$
and so $r(G)=r_0-1$, and then Theorem~\ref{maximum-first-theorem}
implies that $G\iso C(n,r_0-1)$ as desired.

Now that $v_0$ must exist, the possibilities for other blocks in $G$ are severely limited: $G$ can have no other blocks, or a $K_2$, $K_3$, or $K_4$ endblock containing $v_0$, or any number of $D_1$ endblocks which intersect at $v_0$.  Checking the cases where $G$ contains either a complete block or no other block gives us either a graph which is either not extremal or isomorphic to $C(n,r_0-1)$ if that block is isomorphic to $K_3$.  Now suppose that $G$ contains one or more copies of $D_1$ endblocks intersecting at $v_0$.  Let $u_1,u_2,\dots,u_i$ denote the vertices of these blocks that are not adjacent to $v_0$.  Every maximum independent set in $G$ must contain $\{v_0,u_1,u_2,\dots,u_i\}$, and from this it is easy to see that such graphs are not extremal.
This completes the proof of Theorem~\ref{gap2}.
\end{proof}

\bigskip

\noindent{\it Acknowledgment.}  We are indebted to Jason Tedor for
suggesting the use of the Ear Decomposition Theorem, to Herbert Wilf for
suggesting that we look at the maximum independent set problem once we
had done the maximal one, and to the anonymous referees for their numerous
insightful remarks and corrections.

\end{document}